\documentclass{amsart}
\usepackage{amssymb}
\theoremstyle{plain}
\newtheorem{theorem}{Theorem}[section]

\newtheorem{lemma}[theorem]{Lemma}
\newtheorem{proposition}[theorem]{Proposition}

\theoremstyle{definition}
\newtheorem{definition}{Definition}[section]

\theoremstyle{remark}
\newtheorem{remark}{Remark}[section]
\newtheorem{example}{Example}[section]

\numberwithin{equation}{section}

\begin{document}

\author[G.-Y. Chen]{Guan-Yu Chen}
\thanks{The first author is partially supported by NCTS, Taiwan}

\author[L. Saloff-Coste]{Laurent Saloff-Coste}
\thanks{Research partially supported by NSF Grants DMS 0102126, 0306194
and 0603886}

\address{Division of Mathematics, National Center for Theoretical Science, National Tsing Hua
University, Hsinchu 300, Taiwan}
\email{gychen@math.cts.nthu.edu.tw}

\address{Cornell University, Department of Mathematics, Ithaca, NY 14853-4201}
\email{lsc@math.cornell.edu}

\title{The Cutoff phenomenon for randomized riffle shuffles}

\begin{abstract}We study the cutoff phenomenon for generalized riffle shuffles
where, at each step, the deck of cards is cut into a random number of packs
of multinomial sizes which are then riffled together.
\end{abstract}



\maketitle

\section{Introduction} In this article we consider some generalizations
of the standard riffle shuffle of Gilbert, Shannon and Reeds
(GSR-shuffle for short). The GSR-shuffle models  the way
typical card players shuffle cards. First, the deck is cut into two
packs according to an $(n,\frac{1}{2})$-binomial random variable
where $n$ is the number of cards in the deck. Next, cards are
dropped one by one from one or the other pack with probability
proportional to the relative sizes of the packs. Hence, if the
left pack contains $a$ cards and the right pack $b$ cards, the
next card drops from the left pack with probability $a/(a+b)$.

The history of this model is described in \cite[Chap. 4D]{Db}
where the reader will  also find other equivalent definitions and
a discussion of how the model relates to real life card shuffling.
The survey \cite{D-Dur} gives pointers to the many developments
that arose from the study of the GSR model.

Early results concerning the mixing time
(i.e., how many shuffles are needed to mix up the deck) are described in
\cite{A,ADm,Db}. In particular, using ideas of Reeds,
Aldous proved in \cite{A} that, asymptotically as the
number $n$ of cards tends to infinity, it takes $\frac{3}{2}\log_2n$
shuffles to mix up the deck if convergence is measured in total variation
(we use $\log_a$ to denote base $a$ logarithms and $\log$
for natural, i.e., base $e$, logarithms).

In \cite{BD92}, Bayer and Diaconis obtained an exact useful
formula for the probability distribution describing the state of
the deck after $k$ GSR-shuffles. Namely, suppose that cards are
numbered $1$ through $n$ and that we start with the deck in
order. Let $\sigma$ denote a given arrangement of the cards and
let $Q^k_n(\sigma)$ be the probability that the deck is in state
$\sigma$ after $k$ GSR-shuffles. Then
\begin{equation}\label{k-2-shuffles}
Q^k_n(\sigma)= 2^{-kn}{ {n+2^k-r}\choose {n}}\end{equation} where
$r$ is the number of rising sequences in $\sigma$. Given an
arrangement of the deck, a rising sequence is a maximal subset of
cards consisting of successive face values displayed in order. For
instance, the arrangement $3,1,4,5,7,2,8,9,6$ has rising sequences
$(1,2), (3,4,5,6), (7,8,9)$. See \cite{ADm,BD92} for details. By
definition, the total variation distance between two probability
measures $\mu,\nu$ on a set $S$ is given by
\[
    \|\mu-\nu\|_{\mbox{\tiny TV}}=\sup_{A\subset
    S}\{\mu(A)-\nu(A)\}.
\]
Using the formula displayed in (\ref{k-2-shuffles}), Bayer and
Diaconis gave a very sharp version of the fact that the total
variation mixing time is $\frac{3}{2}\log_2n$ for the GSR-shuffle.
\begin{theorem}[Bayer and Diaconis \cite{BD92}] \label{th-BD}
Fix $c\in (-\infty,+\infty)$.
For a deck of $n$ cards, the total variation distance between the
uniform distribution and the distribution of a deck after
$k=\frac{3}{2}\log_2n \;+ c$  {\em GSR}-shuffles is $$
\frac{1}{\sqrt{2\pi}}\int_{-2^{-c}/4\sqrt{3}}^{2^{-c}/4\sqrt{3}}
e^{-t^2/2} dt +O_c(n^{-1/4}).$$
\end{theorem}

This result illustrates beautifully  the so-called cutoff phenomenon
discussed in \cite{A,ADm,AD,Db,Dc,SC,SGO}.
Namely, there is a sharp transition
in convergence to stationarity. Indeed, the integral above
becomes small very fast as $c$ tends to $+\infty$ and gets close to $1$
even faster as $c$ tends to $-\infty$.

The aim of the present paper is to illustrate further the notion
of cutoff using some generalizations of the GSR-shuffle. Along
this way we will observe several phenomena that have not been, to
the best of our knowledge, noticed before. For a deck of $n$ cards
and a given integer $m$, a $m$-riffle shuffle is defined as
follows. Cut the deck into $m$ packs whose sizes $(a_1,\dots,a_m)$
form a multinomial random vector. In other words, the probability
of having packs of sizes $a_1,\dots, a_m$ is
$m^{-n}\frac{n!}{a_1!\dots a_m!}$. Then form a new deck by
dropping cards one by one from these packs with probability
proportional to the relative sizes of the packs. Thus, if the
packs have sizes $(b_1,\dots,b_m)$ then the next card will drop
from pack $i$ with probability $b_i/(b_1+\dots+b_m)$. We will
refer to an $m$-riffle shuffle simply as an $m$-shuffle in what
follows. Obviously the GSR-shuffle is the same as a $2$-shuffle.
A $1$-shuffle leaves the deck unchanged.

These shuffles were considered in
\cite{BD92} where the following two lemmas are proved.
\begin{lemma}\label{L:c2} In distribution,
an $m$-shuffle followed by an independent $m'$-shuffle equals an
$mm'$-shuffle.
\end{lemma}

\begin{lemma}\label{L:c1} For a deck of $n$ cards in order,
the probability that after an $m$-shuffle the deck is in state $\sigma$
depends only of the number $r=r(\sigma)$ of rising sequences of $\sigma$
and equals $Q_{n,m}(r)$ where
\[Q_{n,m}(r)=m^{-n}\binom{n+m-r}{n}.\]
\end{lemma}
For instance, formula (\ref{k-2-shuffles}) for the distribution of
the deck after $k$ GSR-shuffles follows from a direct application
of these two lemmas since $k$ consecutive independent $2$-shuffles
equal a $2^k$-shuffle in distribution. These lemmas will play a
crucial role in this paper as well.

The model we consider  is as follows. Let $p=(p(1),p(2),\dots)$ be
the probability distribution of an integer valued random variable
$X$, i.e., $$P(X=k)=p(k),\;\; k=1,2,\dots.$$ A $p$-shuffle
proceeds by picking an integer $m$ according to $p$ and performing
an $m$-shuffle. In other words, the distribution of a $p$-shuffle
is the $p$-mixture of the $m$-shuffle distributions. Note that
casinos use multiple decks for some games and that these are
shuffled in various ways (including by shuffling machines). The
model above (for some appropriate $p$) is not entirely unrealistic
in this context.

Because of Lemma \ref{L:c1}, the probability  that starting from
a deck in order we obtain a deck
in state $\sigma$ depends only on the
number of rising sequences in $\sigma$ and is given by
\begin{equation}\label{QnX}
Q_{n,p}(r)=\sum_1^\infty p(m)Q_{n,m}(r)=E(Q_{n,X}(r)).
\end{equation}
Abusing notation, if  $\sigma$ denotes a deck arrangement of $n$ cards
with $r$ rising sequences, we write
$$Q_{n,p}(\sigma)=Q_{n,p}(r).$$
Very generally, if $Q$ is a probability measure on deck arrangements
(hence describes a shuffling method), we denote by $Q^k$ the distribution
of the deck after $k$ such shuffles, starting from a deck in order.
For instance, Lemma \ref{L:c2} yields
$$Q_{n,m}^k=Q_{n,m^k}.$$
Let $U_n$ be the uniform distribution on the set of  deck
arrangements of $n$ cards. Although this will not really play a
role in this work, recall that deck arrangements can be viewed as
elements of the symmetric group $S_n$ in such a way that $Q^k$,
the distribution after $k$ successive $Q$-shuffles, is the
$k$-fold convolution of $Q$ by itself. See, e.g.,
\cite{A,BD92,Db,SCfg}. Each of the measures $Q_{n,p}$ generates a
Markov chain on deck arrangements (i.e., on the symmetric group
$S_n$) whose stationary distribution is $U_n$. These chains are
ergodic if $p$ is  not concentrated at $1$. They are not
reversible. Note that \cite{DFP}  studies a similar but different
model based on top $m$ to random shuffles. See \cite[Section
2]{DFP}.

The goal of this paper is to study the convergence of $Q_{n,p}^k$
to the uniform distribution in total variation as $k$ tends to
infinity and, more precisely, the occurrence of a total variation
cutoff for families  of shuffles $\{(S_n,Q_{n,p_n},
U_n)\}_1^\infty$ as the number $n$ of cards grows to infinity and
$p_n$ is a fixed sequence of probability measures on the integers.
To illustrate this, we state the simplest of our results.
\begin{theorem}\label{th-intro}Let $p$ be a probability measure on the positive integers
such that
\begin{equation}\label{def-mu}
\mu=\sum_1^\infty p(k)\log k<\infty.\end{equation} Fix
$\epsilon\in (0,1)$. Then, for any $k_n>
(1+\epsilon)\frac{3}{2\mu}\log n$, we have $$\lim_{n\rightarrow
\infty}\|Q^{k_n}_{n,p}-U_n\|_{\mbox{\em\tiny TV}}=0$$ whereas, for
$k_n<(1-\epsilon)\frac{3}{2\mu}\log n$, $$\lim_{n\rightarrow
\infty}\|Q^{k_n}_{n,p}-U_n\|_{\mbox{\em\tiny TV}}=1.$$
\end{theorem}
In words, this theorem establishes a total variation cutoff at
time $\frac{3}{2\mu}\log n$ (see the definition of cutoff in
Section \ref{sec-cut} below). If $p$ is concentrated at $2$, i.e.,
$Q_{n,p}$ represents a GSR-shuffle, then $\mu=\log 2$ and
$\frac{3}{2\mu}\log n=\frac{3}{2}\log_2n$ in accordance with the
results of Aldous \cite{A} and Bayer-Diaconis \cite{BD92} (e.g.,
Theorem \ref{th-BD}).

The results we obtain are more general and more precise than
Theorem \ref{th-intro} in several directions. First, we will
consider the case where the probability distribution $p=p_n$
depends on the size $n$ of the deck. This is significant because
we will not impose that the sequence $p_n$ converges as $n$ tends
to infinity. Second, and this may be a little surprising at first,
(\ref{def-mu}) is not necessary for the existence of a cutoff and
we will give sufficient conditions that are weaker than
(\ref{def-mu}). Third, under stronger moment assumptions, we will
describe the optimal window size of the cutoff. For instance,
Theorem \ref{th-BD} says that, for the GSR-shuffle, the window
size is of order $1$ with a  normal shape. This result generalizes
easily to any $m$-shuffle where $m$ is a fixed integer greater or
equal to $2$. See Remark \ref{R1} and Theorem \ref{th-BD*} below.
Suppose now that instead of the GSR-shuffle we consider the
$p$-shuffle with $p(2)=p(3)=1/2$. In this case, $\mu=\log
\sqrt{6}$. Theorem \ref{th-intro} gives a total variation cutoff
at time $\frac{3}{2} \log_{\sqrt{6}} n$. We will show that this
cutoff has optimal window size of order $\sqrt{\log n}$. Thus
picking at random between $2$ and $3$ shuffles changes the window
size significantly when compared to either pure $2$-shuffles or
pure $3$-shuffles.

We close this introduction with a remark concerning the spectrum
of these generalized riffle shuffles and how it relates to the
window of the cutoff. As Lemma \ref{L:c2} makes clear, all riffle
shuffles commute. Although riffle shuffles are not reversible,
they are all diagonalizable with real positive eigenvalues and
their spectra can be computed explicitly (this is another
algebraic ``miracle'' attached to these shuffles!). See
\cite{BD92,BHR,BD}. In particular, the second largest eigenvalue
of an $m$-shuffle is $1/m$ with the same eigenspace for all $m\ge
2$. See \cite{M01dr} for a stronger result implying this
statement. Thus, the second largest eigenvalue of a $p$-shuffle is
$\beta=\sum k^{-1}p(k).$ By definition, the relaxation time of a
finite Markov chain is the inverse of the spectral gap
$(1-\beta)^{-1}$ and one might expect that, quite generally, for
families of Markov chains presenting a cutoff, this quantity would
give a good control of the window of the cutoff. The generalized
riffle shuffles studied here provided interesting (albeit
non-reversible) counterexamples: Take, for instance, the case
discussed earlier where $p(2)=p(3)=1/2$. Then $\beta=\frac{5}{12}$
and $(1-\beta)^{-1}= \frac{12}{7}$, independently of the number
$n$ of cards. However, as mentioned above, the optimal window size
of the cutoff for this family is $\sqrt{\log n}$. For generalized
riffle shuffles, the window size of the cutoff and the relaxation
time appear to be disconnected.

\section{The cutoff phenomenon}\label{sec-cut}
The following definition introduces the notion of cutoff for a
family of ergodic Markov chains.
\begin{definition} \label{def-cutoff1}
Let $\{(S_n,K_n,\pi_n)\}_1^\infty$ be a family of ergodic Markov
chains where $S_n$ denotes the state space, $K_n$ the Markov
kernel, and $\pi_n$ the stationary distribution. This family
satisfies a total variation cutoff with critical time $t_n>0$ if,
for any fixed $\epsilon\in (0,1)$, $$\lim_{n\rightarrow
\infty}\sup_{x\in S_n}\|K^{k_n}_n(x,\cdot)-\pi_n \|_{\mbox{\tiny
TV}}=\left\{\begin{array}{cl}0& \mbox{ if }\,
k_n>(1+\epsilon)t_n\\ 1 &\mbox{ if }\,
k_n<(1-\epsilon)t_n.\end{array}\right.$$
\end{definition}
This definition was introduced in \cite{ADm}. A more thorough
discussion is in \cite{Dc} where many examples are described. Note
that this definition does not require that the critical time $t_n$
tends to infinity (in \cite{Dc}, the corresponding definition
requires that $t_n$ tends to infinity). The positive times $t_n$
can be arbitrary and thus can have several limit points
in $[0,\infty]$. Examples of families
having a cutoff with a bounded critical time sequence will be
given below. Theorem \ref{th-intro} above states that, under
assumption (\ref{def-mu}), a $p$-shuffle has a total variation
cutoff with critical time $t_n=\frac{3}{2\mu}\log n$.

Informally, a family has a cutoff if convergence to stationarity
occurs in a time interval of size $o(t_n)$ around the critical
time $t_n$. The size of this time interval can be thought of as
the ``window'' of the cutoff. The next definition carefully
defines the notion of the window size of a cutoff.

\begin{definition}\label{def-cutoff2}
Let $\{(S_n,K_n,\pi_n)\}_1^\infty$ be a family of ergodic Markov
chains as in Definition \ref{def-cutoff1}. We say that this family
presents a $(t_n,b_n)$ total variation cutoff if the following
conditions are satisfied:
\begin{enumerate}
\item For all $n=1,2,\dots,$ we have $t_n>0$ and $\displaystyle
\lim_{n\rightarrow \infty} b_n/t_n=0.$

\item For $c\in\mathbb{R}-\{0\}$ and $n\ge 1$, set
\[
    k=k(n,c)=\begin{cases}\lceil t_n+cb_n\rceil&
    \text{if $c>0$}\\\lfloor
    t_n+cb_n\rfloor&\text{if $c<0$}
    \end{cases}.
\]
The functions $\overline{f},\underline{f}$ defined by
$$\overline{f}(c)=\limsup_{n\rightarrow \infty}\,
\sup_{x\in S_n} \|K^k_n(x,\cdot)-\pi_n\|_{\mbox{\tiny TV}} \text{
for } c\ne 0 $$ and
$$\underline{f}(c)=\liminf_{n\rightarrow \infty}\,
\sup_{x\in S_n} \|K^k_n(x,\cdot)-\pi_n\|_{\mbox{\tiny TV}} \mbox{
for } c\ne 0$$ satisfy
$$\lim_{c\rightarrow \infty} \overline{f}(c)=0,\;\;
\lim_{c\rightarrow -\infty}\underline{f}(c)=1.$$
\end{enumerate}
\end{definition}
\begin{definition}\label{def-opt}
Referring to Definition \ref{def-cutoff2}, a $(t_n,b_n)$ total
variation cutoff is
said to be optimal if the functions $\overline{f},\underline{f}$
satisfy
$\underline{f}(c)>0$ and $ \overline{f}(-c)<1$ for all $c>0$.
\end{definition}
Note that any family having a $(t_n,b_n)$ cutoff
(Definition \ref{def-cutoff2}) has a cutoff with critical time
$t_n$ (Definition \ref{def-cutoff1}). The sequence
$(b_n)_1^\infty$ in Definition \ref{def-cutoff2}
describes an upper bound on the optimal window
size of the cutoff. For instance the main result of Bayer and
Diaconis \cite{BD92}, i.e., Theorem \ref{th-BD} above, shows that
the GSR-shuffle family presents a $(t_n,b_n)$ total variation
cutoff with $t_n=\frac{3}{2}\log_2n$ and $b_n=1$. Theorem
\ref{th-BD} actually determines exactly ``the shape'' of the
cutoff, that is, the two functions $\overline{f},\underline{f}$ of
Definition \ref{def-cutoff2}. Namely, for the GSR-shuffle family
and $t_n=\frac{3}{2}\log_2n$, $b_n=1$, we have $$\overline{f}(c)=
\underline{f}(c)=\frac{1}{\sqrt{2\pi}}\int_{-2^{-c}/4\sqrt{3}}^{2^{-c}/4\sqrt{3}}
e^{-t^2/2} dt.$$ This shows that this cut-off is optimal
(Definition \ref{def-opt}).

The optimality introduced in Definition \ref{def-opt} is very
strong. If a family presents an optimal $(t_n,b_n)$ total
variation cut-off and also a $(s_n,c_n)$ total variation cut-off,
then $t_n\sim s_n$ and $b_n=O(c_n)$. In words, if $(t_n,b_n)$ is
an optimal cut-off then there are no cut-offs with a window
significantly smaller than $b_n$. For a more detailed discussion
of the cutoff phenomena and their optimality, see \cite{C06}.

\section{Cutoffs for generalized riffle shuffles}\label{gen-sh}
In this section we state our main results and illustrate them with
simple examples. They describe total variation cutoffs for
generalized riffle shuffles, that is, for the $p$-shuffles defined
in the introduction. More precisely, for each $n$ ($n$ is the
number of cards), fix a probability distribution
$p_n=(p_n(1),p_n(2),\dots)$ on the integers and consider the
family of Markov chains (i.e., shuffles)
$$\{(S_n,Q_{n,p_n},U_n)\}_1^\infty.$$ Here $S_n$ is the set of all deck
arrangements (i.e., the symmetric group) and $U_n$ is the uniform
measure on $S_n$. For any $x\in
[0,\infty]$, set
\begin{equation}\label{eq-psi}
\Psi(x)=\frac{1}{\sqrt{2\pi}}\int_{-x/4\sqrt{3}}^{x/4\sqrt{3}}
e^{-t^2/2} dt.
\end{equation}

We start with the simple case where
the probability distributions $p_n$
is concentrated on exactly one integer $m_n$ and use the notation $Q_{n,m_n}$
for an $m_n$-shuffle.
\begin{theorem}\label{th-1pt} Let $(m_n)_1^\infty$ be any
sequence of integers all greater than $1$ and set $$\mu_n=\log
m_n,\;\;t_n=\frac{3\log n}{2\mu_n}.$$ Then the family
$\{(S_n,Q_{n,m_n},U_n)\}_1^\infty$ presents a $(t_n, \mu_n^{-1})$
total variation cutoff.
\end{theorem}

\begin{remark} \label{R1}
When $m_n=m$ is constant Theorem \ref{th-1pt} gives a
$(\frac{3}{2}\log_m n,1)$ total variation cutoff. In this case,
for $k=\frac{3}{2}\log_m n \,+ c$, one has the more precise result
that $\|Q^{k}_{n,m}-U_n\|_{\mbox{\tiny TV}}=\Psi(m^{-c})
+O_c(n^{-1/4}).$ In particular, for $m=2$, this is the Theorem of
Bayer and Diaconis stated as Theorem \ref{th-BD} in the
introduction.\end{remark}

Next we give a more explicit version of Theorem
\ref{th-1pt} which requires some additional notation. For any real
$t>0$, set $$\{t\}=\left\{\begin{array}{cl} 1/2 &\mbox{ if } 0<t<
1/2\\
 k &\mbox{ if } k-1/2 \le t < k+1/2 \mbox{ for some  } k=1,2,\dots,
\end{array}\right.$$
(this is a sort of ``integer part'' of $t$) and
$$d(t)=\left\{\begin{array}{cl} 1/2 &\mbox{ if } 0<t< 1/2\\
 t-\{t\} &\mbox{ if } 1/2 \le t <\infty. \end{array}\right.$$

\begin{theorem}\label{th-1pt'} Let $(m_n)_1^\infty$ be any
sequence of integers all greater than $1$. Consider the family of
shuffles $\{(S_n,Q_{n,m_n},U_n)\}_1^\infty$ and let $\mu_n$, $t_n$
be as in Theorem {\em \ref{th-1pt}}.
\begin{itemize}
\item[(A)]
Assume that $\lim\limits_{n\rightarrow\infty}m_n=\infty$, that is,
$\lim\limits_{n\rightarrow\infty}\mu_n=\infty$. Then, we have:
\begin{itemize}
\item[(1)] The family $\{(S_n,Q_{n,m_n},U_n)\}_1^\infty$ always has a
$(\{t_n\},b_n)$ cutoff for any positive $b_n=o(1)$, that is,
 $$\lim_{n\rightarrow \infty}\inf_{k<\{t_n\}}
\|Q^{k}_{n,m_n}-U_n\|_{\mbox{\em \tiny TV}}=1, \;\;
\lim_{n\rightarrow \infty}\sup_{k>\{t_n\}}
\|Q^{k}_{n,m_n}-U_n\|_{\mbox{\em \tiny TV}}=0.$$
\item[(2)]
If $\displaystyle\lim_{n\rightarrow \infty} |d(t_n)|\mu_n=\infty$
then there is a $(t_n,0)$ cutoff, that is,
$$\lim_{n\rightarrow \infty}\inf_{k\le t_n}
\|Q^{k}_{n,m_n}-U_n\|_{\mbox{\em\tiny TV}}=1, \;\;
\lim_{n\rightarrow \infty}\sup_{k\ge t_n}
\|Q^{k}_{n,m_n}-U_n\|_{\mbox{\em \tiny TV}}=0.$$
\item[(3)]
If $\displaystyle\liminf_{n\rightarrow \infty} |d(t_n)|\mu_n
<\infty$ then there exists a sequence $(n_i)_1^\infty$ tending to
infinity such that
\[0< \liminf_{i\rightarrow\infty}
\|Q_{n_i,m_{n_i}}^{\{t_{n_i}\}}
-U_{n_i}\|_{\mbox{\em \tiny TV}}\le
    \limsup_{i\rightarrow\infty}\|Q_{n_i,m_{n_i}}^{\{t_{n_i}\}}-U_{n_i}
\|_{\mbox{\em \tiny TV}}<1.\]
In particular, there is no $(t_n,0)$ total variation
cutoff.
\item[(4)]
If $\displaystyle \lim_{n\rightarrow \infty} d(t_n)\mu_n= L\in
[-\infty,\infty]$ exists then
\begin{equation}\label{eq-th4-4}
    \lim_{n\rightarrow
    \infty}\|Q^{\lfloor\{t_n\}\rfloor}_{n,m_n}-U_n\|_{\mbox{\em \tiny
    TV}}= \Psi(e^{L}).
\end{equation}
\end{itemize}
\item[(B)] Assume that $(m_n)_1^\infty$ is bounded. Then
$t_n$ tends to infinity, there is a $(t_n,1)$ total variation cutoff and,
for any fixed $k\in\mathbb{Z}$, we have
 $$0< \liminf_{n\rightarrow\infty}
\|Q_{n,m_n}^{\{t_n\}+k} -U_{n}\|_{\mbox{\em \tiny TV}}\le
\limsup_{n\rightarrow\infty}\|Q_{n,m_n}^{\{t_n\}+k}-U_n
\|_{\mbox{\em \tiny TV}}<1.$$ In particular, the
$(t_n,1)$ cutoff is optimal.
\end{itemize}
\end{theorem}

\begin{example}\label{E1}
To illustrate this result, consider the case where
$m_n=\lfloor n^\alpha\rfloor
$ for some fixed $\alpha>0$. In this case, we have
$$\mu_n\sim \alpha \log n,\;\;
t_n= \frac{3\log n}{2\mu_n}\sim \frac{3}{2\alpha}
\mbox{ as $n$ tends to infinity}.$$
\begin{itemize}
\item[(a)] Assume that
$\frac{3}{2\alpha}\in (k,k+1)$ for some $k=0,1,2,\dots$.
Then  $|d(t_n)| \mu_n\rightarrow \infty$ and
$$\lim_{n\rightarrow \infty}\|Q^{k}_{n,m_n}-U_n\|_{\mbox{\tiny TV}}=1,
\;\;
\lim_{n\rightarrow \infty}\|Q^{k+1}_{n,m_n}-U_n\|_{\mbox{\tiny TV}}=0.$$
\item[(b)] Assume that  $\frac{3}{2\alpha}=k$ for some integer $k=1,2,\dots$.
Then $|d(t_n)|= O(n^{-\alpha})$. Hence
$|d(t_n)|\mu_n\rightarrow 0$
as $n$ tends to infinity. Theorem \ref{th-1pt'}(1) shows that we have
a $(k, b_n)$ cutoff where $b_n$ is an arbitrary sequence
of positive numbers tending to $0$. That means that
$$\lim_{n\rightarrow \infty}\|Q^{k-1}_{n,m_n}-U_n\|_{\mbox{\tiny TV}}=1,
\;\;
\lim_{n\rightarrow \infty}\|Q^{k+1}_{n,m_n}-U_n\|_{\mbox{\tiny TV}}=0.$$
Moreover Theorem \ref{th-1pt'}(4) gives
$\lim_{n\rightarrow \infty}\|Q^{k}_{n,m_n}-U_n\|_{\mbox{\tiny TV}}=
\Psi(1).$
\end{itemize}
\end{example}

\begin{example}\label{E1'}
Consider the case where $m_n=\lfloor (\log n)^\alpha \rfloor$, $\alpha>0$.
Then
$$\mu_n\sim \alpha \log \log n,\;\;
t_n\sim  \frac{3\log n}{2\alpha \log \log n}
\mbox{ as $n$ tends to infinity}.$$
Note that $t_n$ tends to infinity and the window size $\mu_n^{-1}$
goes to zero.
\end{example}

We now state results concerning general $p$-shuffles.
We will need the following notation. For each $n$, let  $p_n$ be
a probability distribution  on the integers. Let $X_n$
be a random variable with distribution $p_n$. Assume that $p_n$
is not supported on a single integer and set
$$\mu_n=E(\log X_n),\;\;\sigma_n^2=\mbox{Var}(\log X_n),\;\;\xi_n=
\frac{\log X_n \,-\mu_n}{\sigma_n}.$$
Consider the following conditions which may or may not be satisfied by $p_n$:
\begin{equation}\label{mu}
\lim_{n\rightarrow \infty}\frac{\log n}{\mu_n}=\infty.\end{equation}
\begin{equation}\label{Lind} \forall\,\epsilon>0,\;\;\;
\lim_{n\rightarrow \infty} E\left(\xi_n^2
\mathbf{1}_{\{\xi_n^2>\epsilon \mu_n^{-1}\log n\}}\right)=0.
\end{equation}
Condition (\ref{Lind}) should be understood as a Lindeberg type
condition. We will prove in Lemma \ref{L:main3} that (\ref{Lind})
implies (\ref{mu}). 
Example \ref{E4} shows that the converse is false.

\begin{theorem} \label{T:main2}
Referring to the notation introduced above, assume that
$$0<\mu_n,\sigma_n<\infty$$ and set $$ t_n=\frac{3\log
n}{2\mu_n},\;\;
b_n=\frac{1}{\mu_n}\max\left\{1,\sqrt{\frac{\sigma_n^2\log
n}{\mu_n}} \right\}.$$ Assume that the sequence $(p_n)$ satisfies
{\em (\ref{Lind})}. Then the family
$\{(S_n,Q_{n,p_n},U_n)\}_1^\infty$ presents a $(t_n,b_n)$ total
variation cutoff. Moreover, if the window size $b_n$ is bounded
from below by a positive real number, then the $(t_n,b_n)$ total
variation cut-off is optimal.
\end{theorem}

\begin{example}\label{E2}
Assume $p_n=p$ is independent of $n$ and $$\mu=\sum_1^\infty p(k)
\log k <\infty,\;\; \sigma^2 =\sum_1^\infty |\mu-\log k|^2p(k)
<\infty.$$ Then condition (\ref{Lind}) holds and
$$t_n=\frac{3}{2\mu}\log n,\;\; b_n\approx \sqrt{\log n}$$ where
$b_n\approx \sqrt{\log n}$ means that the ratio $b_n/\log n$ is
bounded above and below by positive constants. Thus Theorem
\ref{T:main2} yields an optimal  $(\frac{3}{2\mu}\log n,\sqrt{\log n})$
total variation cutoff.
\end{example}

\begin{example} \label{E3}
Assume that $p_n$ is concentrated equally
on two integers $m_n <m'_n$ and write $m'_n=m_nk^2_n$. Thus
$p_n(m_n)=p_n(m_nk^2_n)=1/2$ and
$$\mu_n=\log m_nk_n, \;\;\sigma_n=\log k_n.$$
In this case, Condition
(\ref{Lind}) is equivalent to (\ref{mu}), that is
$$\mu_n=\log (m_nk_n)= o(\log n).$$
Assuming that (\ref{mu}) holds true, Theorem \ref{T:main2} yields a total
variation cutoff at time
$$t_n=\frac{3\log n}{2\log m_nk_n}$$
with window size
$$b_n=\frac{1}{\log m_nk_n} \max\left\{1,
\sqrt{\frac{(\log k_n)^2\log n}{\log m_n k_n}}\right\}.$$
For instance, assume that $m'_n=m_n+1$ with $m_n$ tending to infinity.
Then (\ref{mu}) becomes $\log m_n=o(\log n)$ and we have
$$b_n= \frac{1}{\log m_n}\max\left\{1,
\frac{(\log n)^{1/2}}{m_n(\log m_n)^{1/2}}\right\}.$$
Specializing further to  
$m_n\approx (\log n)^\alpha$ with $\alpha\in (0,\infty)$ yields
$$t_n\sim \frac{3\log n}{2\alpha \log\log n}$$
and
$$b_n\approx \left\{\begin{array}{cl} (\log\log n)^{-1}
&\mbox{ if } \alpha\in [1/2,\infty)\\
(\log n)^{1/2-\alpha}(\log \log n)^{-3/2} &\mbox{ if } \alpha\in (0,1/2).
\end{array}\right.$$
In particular, $b_n=o(1)$ when $\alpha\ge 1/2$ but tends to
infinity when $\alpha\in (0,1/2)$. Compare with Example \ref{E1'}
above.
\end{example}

Regarding Theorem \ref{T:main2}, one might want to remove the
hypothesis of existence of a second moment  concerning the  random
variables $\log X_n$. It turns out that it is indeed possible but
at the price of losing control of the window of the cutoff. What
may be more surprising is that one can also obtain results without
assuming that the first moment  $\mu_n$ is finite. In some cases,
it might be possible to control the window size by using
convergence to symmetric stable law of exponent $\alpha\in (1,2)$
but we did not pursue this here.

\begin{theorem}\label{T:main3}
Referring to the notation introduced above, assume that $\mu_n>0$
(including possibly $\mu_n=\infty$). Assume further that there
exists a sequence $a_n$ tending to infinity and satisfying
\begin{equation}\label{E:maint301}
a_n=O(\log n),\quad\lim_{n\rightarrow\infty}\frac{(\log n)EZ_n^2}{a_n^2EY_n}=0,
    \quad \lim_{n\rightarrow\infty}\frac{\log
    n}{EY_n}=\infty,
\end{equation}
where $Y_n=Z_n=\log X_n$ if $\log X_n\le a_n$, and $Y_n=0$,
$Z_n=a_n$ if $\log X_n>a_n$. Then the family
$\{(S_n,Q_{n,p_n},U_n)\}_1^\infty$ presents a total variation
cutoff with critical time
$$t_n=\frac{3\log n}{2EY_n}.$$
\end{theorem}

\begin{remark}\label{R4}
In Theorem \ref{T:main3}, if (\ref{E:maint301}) holds for some
sequence $(a_n)$ then it also holds  for any sequence $(d a_n)$
with $d>0$. Moreover, for all $d>0$, $$E\left((\log
X_n)\mathbf{1}_{\{\log X_n\le da_n\}}\right)\sim EY_n.$$ This is
proved in Lemma \ref{equiv} below.
\end{remark}

\begin{example}\label{E4}
Assume $p_n(\lfloor e^i\rfloor)=c_n^{-1}i^{-2}$ for all $1\le i\le
\lfloor\log n \rfloor$, where
$c_n=1+2^{-2}+3^{-2}+\cdots+(\lfloor\log n\rfloor)^{-2}$. Note
that $c_n\rightarrow c=\pi^2/6$ as $n\rightarrow\infty$. In this
case, $\mu_n\sim c^{-1}\log\log n$, $\sigma_n^2\sim c^{-1}\log n$
and for $\epsilon>0$
\[
    E\left[\xi_n^2\mathbf{1}_{\{\xi_n^2<\epsilon\mu_n^{-1}\log
    n\}}\right]\sim \sqrt{\frac{\epsilon}{\log\log n}}.
\]
Hence the Lindeberg type condition (\ref{Lind}) does not hold and
Theorem \ref{T:main2} does not apply. However, if we consider
$a_n=\log n$ and try to apply Theorem \ref{T:main3}, we have
$EY_n=\mu_n\sim c^{-1}\log\log n$ and $EZ_n^2\sim c^{-1}\log n$.
This implies that (\ref{E:maint301}) holds and yields a total
variation cutoff with critical time $\frac{\pi^2\log n}{4\log\log
n}$.

The untruncated version of this example is $p_n(\lfloor
e^i\rfloor)=p(\lfloor e^i\rfloor)=c^{-1}i^{-2}$, $i=1,2,\dots$ and
$c=\pi^2/6$. In this case, $\mu_n=\mu=\infty$. Theorem
\ref{T:main3} applies with $a_n=\log n$ and yields a total
variation cutoff with critical time $\frac{\pi^2\log n}{4\log\log
n}$.\end{example}

We end this section with a result which is a simple corollary of
Theorem \ref{T:main3} and readily implies Theorem \ref{th-intro}.

\begin{theorem} \label{th-intro+}
Let $X_n,p_n,\mu_n$ be as above. Assume that
\begin{equation}\label{hyp1}
\mu_n=E(\log X_n)=o(\log n)
\end{equation}
and that, for any fixed $\eta>0$,
\begin{equation}\label{hyp2}
E[(\log X_n)\mathbf{1}_{\{\log X_n>\eta \log n\}}]=o_\eta(\mu_n).
\end{equation}
Then the family  $\{(S_n,Q_{n,p_n},U_n)\}_0^\infty$ has a total
variation cutoff at time $t_n=\frac{3\log n}{2\mu_n}.$
\end{theorem}

\begin{example} Suppose $p_n=p$ and $0<\mu_n=\mu<\infty$ as in Theorem
\ref{th-intro}. Then condition (\ref{hyp1})-(\ref{hyp2}) are obviously
satisfied. Thus Theorem \ref{th-intro} follows immediately from Theorem
\ref{th-intro+} as mentioned above.
\end{example}

\begin{remark} Condition (\ref{hyp2}) holds true if $X_n$ satisfies
the (logarithmic) moment condition that there exists $\epsilon>0$ such that
$$\frac{E([\log X_n]^{1+\epsilon})}{(\log n)^\epsilon}=o(\mu_n).$$
\end{remark}

\section{An application: Continuous-time card shuffling}
In this section, we consider the continuous-time version of the previous
card shuffling models
where the waiting times between two successive  shuffles are
independent exponential(1) random variables. Thus, the distribution of
card arrangements at time $t$ starting from the deck in
order is given by the
probability measure
$H_{n,t}=e^{-t(I-Q_{n,p_n})}$ defined by
\begin{equation}\label{def-cts}
    H_{n,t}(\sigma)=H_{n,t}(r)=e^{-t}\sum_{k=0}^\infty
    \frac{t^k}{k!} Q_{n,p_n}^k(r)\quad \text{for $\sigma\in S_n$},
\end{equation}
where $r$ is the number of rising sequences of $\sigma$.

The definition of total variation cutoff and its optimality for
continuous time families is the same as in Definitions
\ref{def-cutoff1}, \ref{def-cutoff2} and \ref{def-opt} except that
all times are now taken to be non-negative reals. To state our
results concerning the family $\{(S_n,H_{n,t},U_n)\}_1^\infty$ of
continuous time Markov chains associated with $p_n$-shuffles,
$n=1,2,\dots,$ we keep the notation introduced in Section
\ref{gen-sh}. In particular, we set
$$\mu_n=E(\log X_n),\;\;\sigma_n^2=\mbox{Var}(\log X_n),
\;\;t_n=\frac{3\log n}{2\mu_n},$$
where $X_n$ denotes a random variable with distribution $p_n$,
and, if $\mu_n,\sigma_n\in ( 0,\infty)$,
$$\xi_n=
\frac{\log X_n \,-\mu_n}{\sigma_n}.$$
We will  obtain the following theorems as
corollaries of the discrete time results of  Section
\ref{gen-sh}.
Our first result concerns  the case where each $p_n$ is concentrated
on one integer as in Theorem \ref{th-1pt}.

\begin{theorem}\label{T:mainc1}
Assume that for each $n$ there is an integer $m_n$ such that
$p(m_n)=1$ Then $\mu_n=\log m_n$, $t_n=\frac{3\log n}{2\log m_n}$
and the family $\mathcal{F}=\{(S_n,H_{n,t},U_n)\}_1^\infty$
presents a total variation cutoff if and only if
$$\lim_{n\rightarrow\infty}\frac{\log n}{\log m_n}=\infty.$$
Moreover, if this condition is satisfied then $\mathcal{F}$ has an
optimal $\left(t_n,\sqrt{t_n}\right)$ total variation cutoff.
\end{theorem}

Compare with the discrete time result stated in Theorem \ref{th-1pt}
and with Example \ref{E1} which we now revisit.

\begin{example}
Assume that $P(X_n=\lfloor n^\alpha\rfloor)=1$ for a fixed
$\alpha>0$ as in Example \ref{E1}. According to Theorem
\ref{T:mainc1}, the continuous time family $\mathcal{F}$ does not
present  a total variation cutoff in this case since
$\lim_{n\rightarrow \infty}\frac{\log n}{\mu_n}=1/\alpha<\infty$.
Recall from Example \ref{E1} that the corresponding discrete time
family has a cutoff.

Assume that $P(X_n=
 \lfloor(\log n)^\alpha\rfloor )=1$
for some
fixed $\alpha>0$ as in Example \ref{E2}.
In this case, the family $\mathcal{F}$ presents a
$(t_n,\sqrt{t_n})$ total variation cutoff with
$t_n=\frac{3\log n}{2\alpha\log\log n}.$
Note that the window
of the continuous time cutoff differs greatly from the window of the
discrete time cutoff in this case.
\end{example}

Next we consider the general case under various hypotheses
paralleling Theorems \ref{T:main2} and \ref{T:main3}.
\begin{theorem}\label{T:mainc}
Consider the continuous time family
$\mathcal{F}=\{(S_n,H_{n,t},U_n)\}_1^\infty$ associated to a
sequence $(X_n)_1^\infty$ of integer valued random variables with
probability distributions $(p_n)_1^\infty$.
\begin{itemize}
\item[(1)] Assume that $\mu_n,\sigma_n\in (0,\infty)$ for all $n\ge 1$
and that
{\em (\ref{Lind})} holds. Then the family
$\mathcal{F}$ presents an optimal  $\left(t_n,b_n\right)$
total variation cutoff, where
$$t_n=\frac{3\log n}{2\mu_n},\,\,b_n=\frac{1}{\mu_n}\max\left\{\left(\mu_n+\sigma_n\right)
\sqrt{\frac{\log n}{\mu_n}},1\right\}.$$
\item[(2)] Assume that $\mu_n>0$ (including possibly
$\mu_n=\infty$) and there exists a sequence $(a_n)_1^\infty$
tending to infinity such that {\em (\ref{E:maint301})} holds. Then
$\mathcal{F}$ presents a total variation cutoff with critical time
$$t_n=\frac{3\log n}{2EY_n}$$
where $Y_n=(\log X_n)\mathbf{1}_{\{\log X_n\le a_n\}}$.
\end{itemize}
\end{theorem}

\begin{remark} Theorem \ref{T:mainc}(2) applies  when $p_n=p$
is independent of $n$ and
$\mu=\sum_1^\infty p(k)\log k<\infty$.
In this case,
the family $\mathcal{F}=\{(S_n,H_{n,t},U_n)\}_1^\infty$ presents a total
variation cutoff with critical time $t_n=\frac{3\log n}{2\mu}$ as in Theorem
\ref{th-intro}.  If in addition we assume that
$\sigma^2=\sum_1^\infty|\mu-\log k|^2p(k)<\infty$ then Theorem
\ref{T:mainc}(1) applies and shows that $\mathcal{F}$ has a
$(t_n,\sqrt{\log n})$ total variation cutoff. Compare with Example \ref{E2}.
\end{remark}

We now describe how Theorem \ref{T:mainc} applies to
Examples \ref{E3}-\ref{E4} of Section \ref{gen-sh}.

\begin{example}
Assume, as in Example \ref{E3}, that $p_n(m_n)=p_n(m_nk_n^2)=1/2$.
Assume further that $\mu_n=\log(m_nk_n)=o(\log n)$. Then, by
Theorem \ref{T:mainc}(1), $\mathcal{F}$ presents a
$(t_n,\sqrt{t_n})$ total variation cutoff, where $$t_n=\frac{3\log
n}{2\log m_nk_n}.$$

Finally, for Example \ref{E4}, both in truncated and untruncated
cases, Theorem \ref{T:mainc}(2) implies that the family presents a
total variation cutoff with critical time $\frac{\pi^2\log
n}{4\log\log n}$. However, Theorem \ref{T:mainc}(1) is not
applicable here since, in either case, the Lindeberg type
condition (\ref{Lind}) has been shown failed in Example \ref{E4}.
\end{example}

\section{Technical tools}

Two of the main technical tools we will use have already been stated
as Lemma \ref{L:c2} and \ref{L:c1} in the introduction. In particular,
Lemma \ref{L:c1} gives the probability distribution describing
a deck of $n$ cards after an $m$-shuffle, namely,
$$Q_{n,m}(r)=m^{-n} {{n+m-r}\choose{n}}$$
where $r$ is the number of rising sequences in the arrangement of the deck.
The next three known lemmas give further useful information concerning this
distribution.
\begin{lemma}[Tanny, \cite{T73}] \label{L:c1*}
Let $R_{n,h}$ be the number of deck arrangements
of $n$ cards having $r=n/2+ h$ rising sequences, $1\le r\le n$. Then,
uniformly in $h$,
\[\frac{R_{n,h}}{n!}=\frac{e^{-6h^2/n}}{\sqrt{\pi n/6}}
\left(1+o\left(\frac{1}{\sqrt{n}}\right)\right)\]
\end{lemma}

\begin{lemma}[Bayer and Diaconis, {\cite[Proposition 1]{BD92}}]\label{L:c3}
Fix $a\in (0,\infty)$.
For any integers $n,m$ such that  $c=c(n,m)=mn^{-3/2}>a$ and any
$r=\frac{n}{2}+h\in \{1,2,\dots, n\}$, we have
\begin{align}
    Q_{n,m}\left(\frac{n}{2}+h\right)
=\frac{1}{n!}\exp\bigg\{&\frac{1}{c\sqrt{n}}\bigg(-h+\frac{1}{2}
    +O_a\left(\frac{h}{n}\right)\bigg)\notag\\
    &-\frac{1}{24c^2}-\frac{1}{2}\left(\frac{h}{cn}\right)^2+
O_{a}\left(\frac{1}{c n}\right)\bigg\}\notag
\end{align}
as $n$ goes to infinity.
\end{lemma}

\begin{lemma}[Bayer and Diaconis, {\cite[Proposition 2]{BD92}}]\label{L:c4}
Let $h^*$ be the unique integer
such that $Q_{n,m}\left(\frac{n}{2}+h\right)\ge
\frac{1}{n!}$ if and only if $h\le h^*$.
Fix  $a\in (0,\infty)$.
For any integers $n,m$ such that  $c=c(n,m)=mn^{-3/2}>a$, we have
\[
    h^*=\frac{-\sqrt{n}}{24c}+ O_a\left(1\right)\quad
\]
as $n$ tends to $\infty$.
\end{lemma}
The statements of Lemmas \ref{L:c3} and \ref{L:c4} are somewhat
different from the statement in Propositions 1 and 2 in
\cite{BD92} but the same proofs apply. The following theorem
generalizes \cite[Theorem 4]{BD92}, that is, Theorem \ref{th-BD}
of the introduction. The proof, based on the three lemmas
above, is the same as in \cite{BD92}. It is omitted.
\begin{theorem}\label{th-BD*} Fix $a\in (0,\infty)$.
For any integers $n,m$ such that  $c=c(n,m)=mn^{-3/2}>a$ we have
$$\|Q_{n,m}-U_n\|_{\mbox{\em\tiny TV}}
=\frac{1}{\sqrt{2\pi}}\int_{-1/(4\sqrt{3}\, c)}^{1/(4\sqrt{3}\, c)}
e^{-t^2/2}dt+ O_a\left(n^{-1/4}\right).$$
\end{theorem}

Theorem \ref{th-BD*} provides sufficient information to obtain
good upper bounds on the cutoff times of generalized riffle
shuffles. It is however not sufficient to obtain matching lower
bounds and study the cutoff phenomenon. The reminder of this
section is devoted to results that will play a crucial role in
obtaining sharp lower bounds on cutoff times for generalized
riffle shuffles. It is reasonable to guess that shuffling cards
with an $(m+1)$-shuffle is more efficient than  shuffling cards
with an $m$-shuffle . The following Proposition which is crucial
for our purpose says that this intuition is correct when
convergence to stationarity is measured in total variation.

\begin{proposition}\label{P:int}
For any integers $n,m$, we have
$$\left\|Q_{n,m+1}-U_n\right\|_{\mbox{\em \tiny
TV}}\le\left\|Q_{n,m}-U_n\right\|_{\mbox{\em \tiny TV}}.$$
\end{proposition}

\begin{proof} Let $A_m=\{\sigma\in S_n|Q_{n,m}(\sigma)<
\frac{1}{n!}\}$ for $m\ge 1$. By Lemma \ref{L:main} below, we have
$A_{m+1}\subset A_m$ and $Q_{n,m}(\sigma)\le Q_{n,m+1}(\sigma)$
for $\sigma\in A_{m+1}$. This implies
\begin{align}
    \left\|Q_{n,m}-U_n\right\|_{\mbox{\tiny TV}}
    &=U_n(A_m)-Q_{n,m}(A_m)\ge
    U_n(A_{m+1})-Q_{n,m}(A_{m+1})\notag\\&\ge
    U_n(A_{m+1})-Q_{n,m+1}(A_{m+1})=\left\|Q_{n,m+1}-U_n\right\|_
    {\mbox{\tiny TV}}.\notag
\end{align}
\end{proof}

\begin{lemma}\label{L:main}
For any integers, $n,m$ and $r\in \{1,\dots,n\}$, we have:
\begin{itemize}
\item[(1)] $Q_{n,m}(r)\le Q_{n,m+1}(r)$, if $Q_{n,m}(r)\le
\frac{1}{n!}$.
\item[(2)] $Q_{n,k}(r)>\frac{1}{n!}$ for all $k\ge m$, if
$Q_{n,m}(r)>\frac{1}{n!}$.
\end{itemize}
In particular, if $n,m,r$ are such that $Q_{n,m}(r)\le \frac{1}{n!}$,
then
$$k\mapsto Q_{n,k}(r)$$
is non-decreasing on  $\{1,\dots,m\}$.
\end{lemma}

\begin{proof} We prove this lemma by fixing $n$ and $1\le r\le n$,
and considering all possible cases of $m$. For $1\le m<r$, the
first claim holds immediately from Lemma \ref{L:c1} since
$Q_{n,m}(r)=0$, and no $Q_{n,m}(r)$ satisfies the assumption of
the second claim.

For $m\ge r$, consider the following map
$$x\overset{f}{\longmapsto} n\log\left(\frac{x+1}{x}\right)
+\log\left(\frac{x-r+1}{x-r+1+n}\right)\,\, \forall
x\in[r,\infty).$$ The formula of the distribution of deck
arrangements in Lemma \ref{L:c1} implies
$$f(m)=\log\left(\frac{Q_{n,m}(r)}{Q_{n,m+1}(r)}\right).$$ A
direct computation on the derivative of $f$ shows that
\[
    f'(x)=\frac{n[(2r-n-1)x-(r-1)(r-1-n)]}{x(x+1)(x-r+1)(x-r+1+n)}.
\]

Here we consider all possible relation between $r$ and $n$. If
$r,n$ satisfy $\frac{n+1}{2}\le r\le n$, then the derivative $f'$
is positive on $[r,\infty)$. This implies that $f(x)$ is strictly
increasing for $x\ge r$. As
\begin{equation}\label{eq-lim}
    \lim_{x\rightarrow\infty}f(x)=0,
\end{equation}
it follows that
the function $f$ is negative for $x\ge r$ and hence $Q_{n,m}(r)\le
Q_{n,m+1}(r)$ for $m\ge r$. This proves the first claim. Moreover, as
\begin{equation}\label{eq-lim2}
    \lim_{m\rightarrow\infty}Q_{n,m}(r)=\frac{1}{n!},
\end{equation}
we have $Q_{n,m}(r)\le \frac{1}{n!}$ for all $m\ge r$ and
$\frac{n+1}{2}\le r\le n$.

If $r,n$ satisfy $1\le r< \frac{n+1}{2}$, let
$x_0=\frac{(r-1)(r-1-n)}{2r-n-1}$. In this case, the derivative
$f'$ satisfies
\[
    f'(x)\begin{cases}\ge 0&\text{if $r\le x\le x_0$}\\<0&\text{if
    $x>x_0$}\end{cases}.
\]
This implies that $f$ is either decreasing on $[r,\infty)$ or
increasing on $[r,x_0]$ and decreasing on $(x_0,\infty)$ according
to whether $x_0<r$ or $x_0\ge r$.

On one hand, if $x_0<r$, that is, $f$ is decreasing on
$[r,\infty)$, then (\ref{eq-lim}) implies that $f$ is positive on
$[r,\infty)$, which means, in particular, that $Q_{n,m}(r)\ge
Q_{n,m+1}(r)$ for $m\ge r$. In this case, (\ref{eq-lim2}) implies
that $Q_{n,m}(r)\ge \frac{1}{n!}$ for $m\ge r$.

On the other hand, if $x_0\ge r$, that is, $f$ increases on
$[r,x_0)$ and decreases on $[x_0,\infty)$, then (\ref{eq-lim})
implies that $f$ has at most one zero in $[r,\infty)$.
If $f$ has no zero, then $f$ is positive on $[r,\infty)$
and thus (by (\ref{eq-lim2}))
\[
    Q_{n,m}(r)\ge Q_{n,m+1}(r)\ge \frac{1}{n!}\quad\forall m\ge r,
\]
This proves claim (2) (claim (1) is empty in this case).

If $f$ has a zero, say $z$, then (\ref{eq-lim}) implies that $f<0$
on $[r,z)$ and $f>0$ on $(z,\infty)$. By writing $z=\lfloor
z\rfloor+\epsilon$ with $\epsilon\in[0,1)$, it is easy to check
that
\[
    Q_{n,\lfloor z\rfloor}(r)=Q_{n,\lfloor
    z\rfloor+1}(r)>Q_{n,m}(r),\,\, \forall m\ge r,\,m\notin\{\lfloor
    z\rfloor,\lfloor z\rfloor+1\},
\]
when $\epsilon=0$, and
\[
    Q_{n,\lfloor z\rfloor+1}(r)>Q_{n,m}(r),\, \forall m\ge r,\,m\ne
    \lfloor z\rfloor+1,
\]
when $\epsilon\in(0,1)$. Moreover, if  we set $m_0=\lfloor
z\rfloor +1$, then the map $m\mapsto Q_{n,m}(r)$ is increasing on
$[r,m_0]$ and strictly decreasing on $[m_0,\infty)$. In the region
$[m_0,\infty)$, (\ref{eq-lim2}) implies as before that
$Q_{n,m}(r)>\frac{1}{n!}$ for $m\ge m_0$. In the region $[r,m_0]$,
let $m_1\ge r$ be the largest integer $m$ such that $Q_{n,m}(r)\le
\frac{1}{n!}$. Then the monotonicity of the map $m\mapsto
Q_{n,m}(r)$ implies that $Q_{n,m}(r)\le\frac{1}{n!}$ for $r\le
m\le m_1$ and $Q_{n,m}(r)>\frac{1}{n!}$ for $m_1<m\le m_0$. This
proves the desired inequalities.
\end{proof}

\begin{lemma}\label{L:main2} Consider all deck arrangements of a
deck of $n$ cards.
\begin{itemize}
\item[(1)] For $1\le r\le n$, let $A_r$ be the set of deck arrangements with
at least $r$ rising sequences. Then for all integers $n,m$ and
$r\in\{1,\dots, n\}$, we have
$$U_n(A_r)-Q_{n,m}(A_r)\ge
0.$$
\item[(2)]
Fix $a>0$. For integers $n,m$, let $c=c(n,m)=mn^{-3/2}>a$. Let
$B_c$ be the set of deck arrangements with number of rising
sequences in
$[\frac{n}{2}-\frac{\sqrt{n}}{24c}+n^{\frac{1}{4}},n]$. Then
\[
\inf_{k\le m}\bigg(U_n(B_c)-Q_{n,k}(B_c)\bigg)=
\frac{1}{\sqrt{2\pi}}\int_{-1/(4\sqrt{3}\,c)}^{1/(4\sqrt{3}\,c)}
e^{-t^2/2} dt +O_a\left(n^{-\frac{1}{4}}\right).
\]
\end{itemize}
\end{lemma}

\begin{proof} As $Q_{n,m}(r)$ is
non-increasing in $r$, we have either $Q_{n,m}(\sigma)\le
\frac{1}{n!}$ for all $\sigma\in A_r$ or $Q_{n,m}(\sigma)\ge
\frac{1}{n!}$ for all $\sigma\in S_n-A_r$. The inequality stated in (1)
thus follows from the obvious identity
\[
    U_n(A_r)-Q_{n,m}(A_r)=Q_{n,m}(S_n-A_r)-U_n(S_n-A_r).
\]

To prove (2), let $h_0=-\frac{\sqrt{n}}{24c}+n^{\frac{1}{4}}$. By Lemma
\ref{L:c4}, since $h_0\ge h^*$ for large $n$, we have
$Q_{n,m}(\sigma)\le \frac{1}{n!}$ for $\sigma\in B_c$. Lemma
\ref{L:main} then implies
\[
    \inf_{k\le
    m}\bigg(U_n(B_c)-Q_{n,k}(B_c)\bigg)=U_n(B_c)-Q_{n,m}(B_c)
    \quad\text{for $n$ large}.
\]
By Lemmas \ref{L:c1*}, \ref{L:c4}, we have
\begin{align}
    &\left|\bigg(U_n(B_c)-Q_{n,m}(B_c)\bigg)-\|Q_{n,m}-U_n\|_{\mbox{\tiny
    TV}}\right|\le
    \sum_{h=h^*}^{h_0}\frac{R_{n,h}}{n!}\notag\\
    =&\frac{1}{\sqrt{2\pi}}\int_{h^*\sqrt{12/n}}^{h_0\sqrt{12/n}}
    e^{-t^2/2}dt+O\left(n^{-\frac{1}{2}}\right)
    =O_a\left(n^{-\frac{1}{4}}\right).\notag
\end{align}
The equality in (2) then follows from
Theorem \ref{th-BD*}.
\end{proof}

\section{Proof of Theorem \ref{th-1pt}, \ref{th-1pt'}}

The following lemma is a corollary of Theorem
\ref{th-BD*}. It
is the main tool used to prove Theorems
\ref{th-1pt} and \ref{th-1pt'}.

\begin{lemma}\label{C:int2}
For $n\in\mathbb{N}$, let $m_n\in\mathbb{N}$ and
$c_n=m_nn^{-3/2}$. Set
$$\liminf\limits_{n\rightarrow\infty}c_n=L,\,\,
\limsup\limits_{n\rightarrow\infty}c_n=U.$$
\begin{itemize}
\item[(1)] If $L>0$(including possibly the infinity),
then $$\limsup_{n\rightarrow\infty}\|Q_{n,m_n}-U_n\|
_{\mbox{\em\tiny TV}} \le \Psi(L^{-1}).$$

\item[(2)] If $U<\infty$(including possibly 0), then
$$\liminf_{n\rightarrow\infty} \|Q_{n,m_n}-U_n\|_{\mbox{\em\tiny
TV}}\ge \Psi(U^{-1}).$$

\item[(3)] If $U=L\in[0,\infty]$, then
$$\lim_{n\rightarrow\infty}\|Q_{n,m_n}-U_n\|_{\mbox{\em\tiny TV}}
=\Psi(U^{-1}).$$
\end{itemize}
\end{lemma}

\begin{proof}
Note that (3) follows immediately from (1) and (2).  As
the proofs of (1) and (2) are similar, we only prove (1).
Assume first that $0<L<\infty$. Let $\epsilon\in(0,L)$
and choose $N=N(\epsilon)$ such that $c_n\ge L-\epsilon$ for $n\ge
N$. This implies that for $n\ge N$,
\begin{align}
    \|Q_{n,m_n}-U_n\|_{\mbox{\tiny TV}}&\le \sup_{k\ge
    (L-\epsilon)n^{3/2}}\|Q_{n,k}-U_n\|_{\mbox{\tiny TV}}\notag\\
    &=\Psi((L-\epsilon)^{-1})+O_L\left(n^{-1/4}\right),\notag
\end{align}
where the last equality follows from Theorem \ref{th-BD*}.
Letting $n$ tend to infinity first and then
$\epsilon$  to $0$ gives (1).

If $L=\infty$, let $C\in(0,\infty)$ and choose $N=N(C)$ so large
that $c_n\ge C$ if $n\ge N$. As in the previous case, for $n\ge N$,
$$\|Q_{n.m_n}-U_n\|_{\mbox{\tiny TV}}\le
\Psi(C^{-1})+O_C\left(n^{-1/4}\right).$$ Now letting $n,C$ tend to
infinity yields (1) again.

\end{proof}

\noindent{\bf Proof of Theorem \ref{th-1pt}.} For $n\ge 1$ and
$c\in\mathbb{R}$, let $t_n=\frac{3\log n}{2\mu_n}$ and
\[
    k=k(n,c)=\begin{cases}\lceil t_n+c\mu_n^{-1}\rceil&\text{if
    $c>0$}\\ \lfloor t_n+c\mu_n^{-1}\rfloor&\text{if
    $c<0$}\end{cases}.
\]
This implies
\[
    m_n^{k}n^{-3/2}\begin{cases}\ge e^c&\text{if $c>0$}\\
    \le e^c&\text{if $c<0$}\end{cases}.
\]
Let $\overline{f},\underline{f}$ be the functions introduced in
Definition \ref{def-cutoff2}. By Lemmas \ref{L:c2} and
\ref{C:int2}, we have
\[
    \overline{f}(c)\le \Psi(e^{-c})\quad \text{if $c>0$},
\]
and
\[
    \underline{f}(c)\ge\Psi(e^{-c})\quad\text{if $c<0$}.
\]
Letting $c$ tend respectively to $\infty$ and $-\infty$ proves
Theorem \ref{th-1pt}. \hfill $\Box$\vskip2mm

\noindent{\bf Proof of Theorem \ref{th-1pt'}.} In this proof, $k$
always denotes a non-negative integer. We first assume that $m_n$
tends to infinity. Note that
\[
k\begin{cases}\ge t+1/2&\text{if $k>\{t\}$}\\ \le t-1/2&\text{if
$k<\{t\}$ and $t\in[1/2,\infty)$}\\ =0&\text{if $k<\{t\}$ and
$t\in(0,1/2)$}\end{cases}.
\]
This implies
\[
    m_n^{k}n^{-3/2}\begin{cases}\ge m_n^{1/2}&\text{if
    $k>\{t_n\}$}\\ \le m_n^{-1/2}&\text{if $k<\{t_n\}$ and $t_n\ge 1/2$}\\
    =n^{-3/2}&\text{if $k<\{t_n\}$ and
    $t_n\in(0,1/2)$}\end{cases}.
\]
Theorem \ref{th-1pt'}(1) thus follows from  Lemmas \ref{L:c2} and
\ref{C:int2}.

The proof of Theorem \ref{th-1pt'}(2)
is similar to the proof of (1) but depends on the
observation that
\[
    k\begin{cases}\ge t+|d(t)|&\text{if $k>t$}\\ \le
    t-|d(t)|&\text{if $k<t$}\end{cases}\quad \text{for}\,\, k\in\mathbb{N},
\]
which implies
\[
    m_n^kn^{-3/2}\begin{cases}\ge \exp\{|d(t_n)|\mu_n\}&\text{if $k>t_n$}
    \\ \le \exp\{-|d(t_n)|\mu_n\}&\text{if $k<t_n$}.\end{cases}
\]

For Theorem \ref{th-1pt'}(3), by assumptions
$$\liminf_{n\rightarrow\infty}|d(t_n)|\mu_n<\infty,\,\,
\lim_{n\rightarrow\infty}\mu_n=\infty.$$ Thus we can choose $M>0$ and a
sequence $(n_i)_1^\infty$ tending to infinity such that
$|d(t_{n_i})|\mu_{n_i}\le M$ and $t_{n_i}\ge 1/2$ for all $i\ge
1$. Since $\{t\}=t-d(t)$ for $t\ge 1/2$, we have that for all
$i\ge 1$,
\[
    e^{-M}\le m_{n_i}^{\{t_{n_i}\}}n_i^{-3/2}\le e^M.
\]
By Lemmas \ref{L:c2} and \ref{C:int2}, this implies that
\[
    \limsup_{i\rightarrow\infty}\|Q_{{n_i},m_{n_i}}^{\{t_{n_i}\}}-U_{n_i}\|
    _{\mbox{\tiny TV}}\le \Psi(e^{M})<1,
\]
and
\[
    \liminf_{i\rightarrow\infty}\|Q_{{n_i},m_{n_i}}^{\{t_{n_i}\}}-U_{n_i}\|
    _{\mbox{\tiny TV}}\ge \Psi(e^{-M})>0.
\]

For Theorem \ref{th-1pt'}(4), if $L<\infty$, then the fact,
$\lim_{n\rightarrow\infty}\mu_n=\infty$, implies that $t_n\ge 1/2$
for large $n$. In this case, $\{t_n\}=t_n-d(t_n)\in\mathbb{Z}$ and
\begin{equation}\label{eq-mn}
    m_n^{\{t_n\}}=n^{3/2}e^{-d(t_n)\mu_n}.
\end{equation}
Then the desired inequality (\ref{eq-th4-4}) follows from Lemmas
\ref{L:c2} and \ref{C:int2}.

If $L=\infty$, let $(n_i)_1^\infty$ be a sequence such that
$t_n\ge 1/2$ if and only if $n=n_i$ for some $i$. Observe that if
$n\notin \{n_i|i=1,2,...\}$, then $\lfloor\{t_n\}\rfloor=0$, and
hence (\ref{eq-th4-4}) follows immediately. For the sequence
$(t_{n_i})_1^\infty$, since (\ref{eq-mn}) holds in this case, the
discussion for $L<\infty$ is applicable for $t_{n_i}$ and hence
(\ref{eq-th4-4}) holds. This finishes the proof of (4).

We now assume that $(m_n)_1^\infty$ is bounded and let $N$ be an
upper bound of $m_n$. The proof in this case is similar to the
proof of (3) after observing that
\[
    t_n+k-1<\{t_n\}+k< t_n+k+1,
\]
and
\[
    \min\{2^{k-1},N^{k-1}\}\le m_n^{\{t_n\}+k}n^{-3/2}\le
    \max\{2^{k+1},N^{k+1}\}.
\]
\hfill$\Box$\vskip2mm

\section{Proof of Theorem \ref{T:main2}}

We start with the following elementary but crucial lemma.
\begin{lemma}\label{L:main3}
Let $\{Y_n\}_{n=1}^{\infty}$ be a sequence of nonnegative random
variables. Set
$$\mu_n=E[Y_n], \;\;\sigma_n^2=\text{\normalfont
Var}(Y_n),\;\;\xi_n=\frac{Y_n-\mu_n}{\sigma_n}.$$
Suppose that $(a_n)_{n=1}^{\infty}$ is a sequence of
positive numbers such that the Lindeberg type condition
\begin{equation}\label{L1} \forall\, \epsilon>0,\;
    \lim_{n\rightarrow\infty}E\left[\xi_n^21_{\left\{\xi_n^2>\epsilon
    a_n\right\}}\right]=0,
\end{equation}
 holds. Then
\[
    \lim_{n\rightarrow\infty}a_n=\infty\quad \text{and}\quad \lim_{n\rightarrow\infty}\frac{\sigma_n^2}{\mu_n^2a_n}=0.
\]
\end{lemma}

\begin{proof} Note that
$E\left[\xi_n^21_{\{\xi_n^2\le\epsilon a_n\}}\right]\le \epsilon
a_n$ for all $\epsilon>0$. By (\ref{L1}), this implies
$$\liminf\limits_{n\rightarrow\infty} a_n\ge \epsilon^{-1}E[\xi_n^2]=
\epsilon^{-1}.$$
Hence $\lim_{n\rightarrow \infty}a_n=\infty$.
Next, fix $\epsilon>0$.  As $Y_n$ is
nonnegative, we have
\[
    E\left[\xi_n^21_{\{\xi_n<0\}}\right]\le
    \frac{\mu_n^2}{\sigma_n^2}\le \frac{\sqrt{\epsilon}\mu_n^2a_n}{\sigma_n^2},
\]
for all $n$ large enough, and
\[
    E\left[\xi_n^21_{\{0<\xi_n\le \sqrt{\epsilon
    a_n}\}}\right]\le\sigma_n^{-1}\sqrt{\epsilon
    a_n}E\left[(Y_n-\mu_n)1_{\{0<\xi_n\le \sqrt{\epsilon
    a_n}\}}\right]\le \sqrt{\frac{\epsilon\mu_n^2a_n}{\sigma_n^2}}.
\]
Let
$L=\liminf\limits_{n\rightarrow\infty}\mu_n^2a_n/\sigma_n^2\in[0,\infty]$.
Combining both inequalities and letting $n\rightarrow\infty$ imply
\[
    1\le\sqrt{\epsilon}(L+\sqrt{L}).
\]
Letting $\epsilon\rightarrow 0$ shows that $L=\infty$, that is,
$\sigma_n^2/(\mu_n^2a_n)\rightarrow 0$.
\end{proof}


Recall the generalized model of riffle shuffle defined in
(\ref{QnX}). For $n\ge 1$, let $p_n$ be the distribution of an
integer-valued random variable $X_n$ and consider the family
$\{(S_n,Q_{n,p_n},U_n)\}_1^\infty$ where
\[
    Q_{n,p_n}(\cdot)=E(Q_{n,X_n}(\cdot))
    =\sum_{m=1}^\infty p_n(m)Q_{n,m}(\cdot).
\]
Let $X_{n,1},X_{n,2},...$ be a sequence of i.i.d. random variables
sharing the same distribution as $X_n$. Then, for $a,k>0$,
\begin{equation}\label{eq-upper}
\begin{aligned}
    &\|Q_{n,p_n}^k-U_n\|_{\mbox{\tiny TV}}\le \sum_{m=1}^\infty
    P\left(\prod_{i=1}^kX_{n,i}=m\right)\|Q_{n,m}-U_n\|_{\mbox{\tiny TV}}\\
    \le &P\left(\prod_{i=1}^kX_{n,i}\le n^{3/2}a\right)+
    P\left(\prod_{i=1}^kX_{n,i}\ge
    n^{3/2}a\right)\left(\Psi(a^{-1})+O_a\left(n^{-1/4}\right)\right)\\
    =&(\Psi(a^{-1})-1)P\left\{\prod_{i=1}^kX_{n,i}\ge
    n^{3/2}a\right\}+1+O_a\left(n^{-1/4}\right),
\end{aligned}
\end{equation}
where the first inequality comes from the triangle inequality and
the second inequality follows from Theorem \ref{th-BD*}.

Consider the set $B_a$ defined in Lemma \ref{L:main2}, that is,
the subset of $S_n$ containing permutations with numbers of rising
sequences in $[\frac{n}{2}-\frac{\sqrt{n}}{24a}+n^{1/4},n]$. Lemma
\ref{L:main2} then implies that
\begin{equation}\label{eq-lower}
\begin{aligned}
    &\|Q_{n,p_n}^k-U_n\|_{\mbox{\tiny TV}}\ge \sum_{m\le n^{3/2}a}
    P\left\{\prod_{i=1}^kX_{n,i}=m\right\}(U_n(B_a)-Q_{n,m}(B_a))\\
    \ge& \Psi(a^{-1})P\left\{\prod_{i=1}^kX_{n,i}\le
    n^{3/2}a\right\}+O_a\left(n^{-1/4}\right).
\end{aligned}
\end{equation}

\noindent{\bf Proof of Theorem \ref{T:main2}.} For
$c\in\mathbb{R}-\{0\}$, let
\[
    k=k(n,c)=\begin{cases}
    \lceil t_n+cb_n\rceil
    &\text{if $c>0$}\\
    \lfloor t_n+cb_n\rfloor
    &\text{if $c<0$}\end{cases},
\]
where $t_n=\frac{3\log n}{2\mu_n}$ and
$b_n=\frac{1}{\mu_n}\max\left\{1,\sqrt{\frac{\sigma_n^2\log
n}{\mu_n}}\right\}$. By hypothesis,  (\ref{Lind}) holds.
Thus  Lemma \ref{L:main3} implies
\begin{equation}\label{eq-tnbn}
    \lim_{n\rightarrow\infty}t_n=\infty,\,\,
    b_n=o(t_n).
\end{equation}
By Definition \ref{def-cutoff2},
to prove a $(t_n,b_n)$ total variation cut-off,
we have to show that
\[
    \lim_{c\rightarrow\infty}\overline{f}(c)=0\,\,
    \lim_{c\rightarrow-\infty}\underline{f}(c)=1,
\]
where
\[
    \overline{f}(c)=\limsup_{n\rightarrow\infty}
    \|Q_{n,p_n}^{k}-U_n\|_{\mbox{\tiny TV}},\,\,
    \underline{f}(c)=\liminf_{n\rightarrow\infty}
    \|Q_{n,p_n}^{k}-U_n\|_{\mbox{\tiny TV}}.
\]

Note that $b_n\ge
\frac{1}{2\mu_n}\left(1+\sqrt{\frac{\sigma_n^2\log
n}{\mu_n}}\right)$. This implies
\[
    \log(n^{3/2}e^{c/2})-k\mu_n+
    \frac{c}{2}\sqrt{\frac{\sigma_n^2\log n}{\mu_n}}\begin{cases}
    \le 0&\text{if $c>0$}\\\ge 0&\text{if $c<0$}\end{cases}.
\]
Hence, we have
\[
    P\left\{\prod_{i=1}^{k}X_{n,i}\ge n^{\frac{3}{2}}e^{c/2}\right\}
    \ge P\left\{\frac{\sum_{i=1}^{k}\log X_{n,i} -k\mu_n}
    {\sigma_n\sqrt{k}} \ge
    -\frac{c}{2}\sqrt{\frac{\log n}{k\mu_n}}\right\}\,\,
    \text{for $c>0$}
\]
and
\[
    P\left\{\prod_{i=1}^{k}X_{n,i}\le n^{\frac{3}{2}}e^{c/2}\right\}
    \ge P\left\{\frac{\sum_{i=1}^{k}\log X_{n,i} -k\mu_n}
    {\sigma_n\sqrt{k}} \le
    -\frac{c}{2}\sqrt{\frac{\log n}{k\mu_n}}\right\}\,\,
    \text{for $c<0$.}
\]

For fixed $c\in\mathbb{R}-\{0\}$, consider a triangular array of
random variables whose $k$-th row consists of
$$\log X_{n,1},\log X_{n,2},...,\log X_{n,k}.$$ In this
setting, $k\sim t_n$ and (\ref{Lind}) is equivalent to the
well-known Lindeberg condition for such an array. Hence the
central limit theorem (e.g., \cite[Theorem 1, page 329]{S96})
yields
\[
    \liminf_{n\rightarrow\infty}
    P\left\{\prod_{i=1}^{k}X_{n,i}\ge
    n^{\frac{3}{2}}e^{c/2}\right\}\ge
    \frac{1}{2}\left(1+\Psi(2\sqrt{2}c)\right)\quad\text{if
    $c>0$},
\]
and
\[
    \liminf_{n\rightarrow\infty}
    P\left\{\prod_{i=1}^{k}X_{n,i}\le
    n^{\frac{3}{2}}e^{c/2}\right\}\ge
    \frac{1}{2}\left(1+\Psi(-2\sqrt{2}c)\right)\quad\text{if
    $c<0$}.
\]

Then, by (\ref{eq-upper}) and (\ref{eq-lower}), we have
\[
    \overline{f}(c)\le 1-\frac{1}{2}\left(1-\Psi(e^{-c/2})\right)
    \left(1+\Psi(2\sqrt{2}c)\right)\quad\text{for $c>0$},
\]
and
\[
    \underline{f}(c)\ge \frac{1}{2}\Psi(e^{-c/2})
    \left(1+\Psi(-2\sqrt{2}c)\right)\quad\text{for $c<0$}.
\]
Hence the $(t_n,b_n)$-cutoff is proved by letting $c$ tend to
$\infty$ and $-\infty$ respectively.

For the optimality of such total variation cutoff, we need to
estimate $\overline{f}(c)$ for $c<0$ and $\underline{f}(c)$ for
$c>0$. Assume that $b_n\ge b>0$ for all $n\ge 1$. Then we have
\[
    k=\begin{cases}
    \lfloor t_n+cb_n\rfloor>t_n+cb_n-1\ge
    t_n+(c-b^{-1})b_n&\text{if $c<0$}\\
    \lceil t_n+cb_n\rceil<t_n+cb_n+1\le
    t_n+(c+b^{-1})b_n&\text{if $c>0$}\end{cases}.
\]

Arguing as in the proof of cutoff above, we obtain
 \[
    \liminf_{n\rightarrow\infty}P\left\{\prod_{i=1}^{k}X_{n,i}\ge
    n^{\frac{3}{2}}e^{(c-b^{-1})}\right\}\ge
    \frac{1}{2}\left(1-\Psi(4\sqrt{2}(b^{-1}-c))\right)\,\,
    \text{for $c<0$},
 \]
 and
 \[
    \liminf_{n\rightarrow\infty}P\left\{\prod_{i=1}^{k}X_{n,i}\le
    n^{\frac{3}{2}}e^{(c+b^{-1})}\right\}\ge
    \frac{1}{2}\left(1-\Psi(4\sqrt{2}(b^{-1}+c))\right)\,\,
    \text{for $c>0$}.
 \]

Hence, the functions $\overline{f},\underline{f}$ are bounded by
 \[\forall\,c<0,\;\;
    \overline{f}(c)\le
    1-\frac{1}{2}\left(1-\Psi(e^{(b^{-1}-c)})\right)
    \left(1-\Psi(4\sqrt{2}(b^{-1}-c))\right)<1,
 \]
 and
 \[\forall\,c>0,\;\;
    \underline{f}(c)\ge\frac{1}{2}\Psi(e^{-(b^{-1}+c)})
    \left(1-\Psi(4\sqrt{2}(b^{-1}+c))\right)>0.
 \]

 By Definition \ref{def-opt}, the family
 $\{(S_n,Q_{n,p_n},U_n)\}_1^\infty$ has an optimal $(t_n,b_n)$
  total variation cutoff.

\hfill$\Box$\vskip2mm

\section{Proof of Theorems \ref{T:main3}, \ref{th-intro+}}

To work without assuming the existence of $\mu_n$,
we need the following weak law
of large numbers for triangular arrays. See, e.g., \cite{D95}.

\begin{theorem}(Weak law of large numbers)\label{T:wlln}
For each $n$, let $W_{n,k}$, $1\le k\le n$, be independent. Let
$b_n>0$ with $b_n\rightarrow\infty$, and
$\bar{W}_{n,k}=W_{n,k}1_{\{|W_{n,k}|\le b_n\}}$. Suppose that

\noindent{\em (1)} $\sum_{k=1}^nP\{|W_{n,k}|>b_n\}\rightarrow 0$, and

\noindent{\em (2)} $b_n^{-2}\sum_{k=1}^nE\bar{W}_{n,k}^2\rightarrow 0$
as $n\rightarrow\infty$.

If we set $S_n=W_{n,1}+...+W_{n,n}$ and put
$s_n=\sum_{k=1}^nE\bar{W}_{n,k}$, then
\[
    \frac{S_n-s_n}{b_n}\rightarrow 0\quad\text{in probability}.
\]

\end{theorem}

\noindent{\bf Proof of Theorem \ref{T:main3}.}
For $0<|\epsilon|<1$, let
\[
    k=k(n,\epsilon)=\begin{cases}\lceil(1+\epsilon)t_n\rceil&
    \text{if $\epsilon>0$}\\\lfloor(1+\epsilon)t_n\rfloor&\text{if
    $\epsilon<0$}\end{cases}.
\]
By (\ref{eq-upper}) and (\ref{eq-lower}), to prove a total
variation cutoff with critical time $t_n$, it suffices to prove
that for all $a>0$
\begin{equation}\label{eq-upper*}
    \lim_{n\rightarrow\infty}P\left\{\prod_{i=1}^kX_{n,i}
    \ge n^{\frac{3}{2}}a\right\}=1,\,\,
    \text{if $\epsilon>0$},
\end{equation}
and
\begin{equation}\label{eq-lower*}
    \lim_{n\rightarrow\infty}P\left\{\prod_{i=1}^kX_{n,i}
    \le n^{\frac{3}{2}}a\right\}=1,\,\,
    \text{if $\epsilon<0$}.
\end{equation}
Indeed, if these limits holds true then (\ref{eq-upper}) and (\ref{eq-lower})
give
\[
    \limsup_{n\rightarrow\infty}\|Q_{n,p_n}^k-U_n\|_{\mbox{\tiny
    TV}}\le \Psi(a^{-1})\quad \text{for $\epsilon>0$}
\]
and
\[
    \liminf_{n\rightarrow\infty}\|Q_{n,p_n}^k-U_n\|_{\mbox{\tiny
    TV}}\ge \Psi(a^{-1})\quad \text{for $\epsilon<0$}.
\]
The total variation cutoff is then proved by letting $a$ tend to
infinity and 0 respectively.

To prove (\ref{eq-upper*})-(\ref{eq-lower*}),
note that $EZ_n^2=EY_n^2+a_n^2P\{\log X_n>a_n\}$.
By the second part of assumption (\ref{E:maint301}), we have
\begin{equation}\label{wkl*}
    (1+\epsilon)t_nP\{\log X_n>a_n\}\rightarrow
    0\quad\text{and}\quad
    (1+\epsilon)t_na_n^{-2}EY_n^2\rightarrow 0,\quad \text{as
    $n\rightarrow\infty$}.
\end{equation}
In order to apply Theorem \ref{T:wlln}, for fixed
$\epsilon\in(-1,1)$, consider
$$W_{k,1}=\log X_{n,1},...,W_{k,k}=\log X_{n,k}$$
as the $k$-th row of a triangular array of random variables. Then
(\ref{wkl*}) shows that the hypotheses (1) and (2) in Theorem
\ref{T:wlln} hold. Hence
\begin{equation}\label{eq-inprob}
    a_n^{-1}\left(\sum_{i=1}^k\log
    X_{n,i}-(1+\epsilon)t_nEY_n\right)\rightarrow 0\,\,
    \text{in probability}.
\end{equation}
Note also that for $a>0$,
$a_n^{-1}\left(\log(n^{3/2}a)-(1+\epsilon)t_nEY_n\right)\sim
\frac{-3\epsilon\log n}{2a_n}$. Hence the first part of assumption
(\ref{E:maint301}) implies that
\begin{equation}\label{eq-limsupinf}
\begin{aligned}
    &\limsup_{n\rightarrow\infty}a_n^{-1}\left(\log(n^{3/2}a)
    -(1+\epsilon)t_nEY_n\right)<0\quad
    \text{if $\epsilon>0$},\\
    &\liminf_{n\rightarrow\infty}a_n^{-1}\left(\log(n^{3/2}a)
    -(1+\epsilon)t_nEY_n\right)>0\quad
    \text{if $\epsilon<0$}.
\end{aligned}
\end{equation} Combining both (\ref{eq-inprob}) and
(\ref{eq-limsupinf}) proves (\ref{eq-upper*}) and
(\ref{eq-lower*}).\hfill$\Box$\vskip2mm

\noindent{\bf Proof of Theorem \ref{th-intro+}.} Let $X_n$ be
integer valued random variables such that
\[
    P\{X_n=k\}=p_n(k)\quad \text{for } k=1,2,...
\]
and satisfying (\ref{hyp1}), (\ref{hyp2}).
Let $a_n=\log n$ in Theorem \ref{T:main3} so that
$$Y_n=(\log X)\mathbf{1}_{\{\log X\le\log n\}},\,\, Z_n=Y_n+(\log
n)\mathbf{1}_{\{\log X>\log n\}}.$$
Set $L_n=\log X_n$. By (\ref{hyp2}), we have
$E(L_n\mathbf{1}_{\{L_n>\log n\}})=o(\mu_n).$
Hence $E(Y_n)\sim \mu_n$ and the third condition of
(\ref{E:maint301}) follows from (\ref{hyp1}).
To apply Theorem \ref{T:main3}, it remains to show
\[
    \lim_{n\rightarrow\infty}\frac{EY_n^2}{EY_n\log n}=0,\quad
    \lim_{n\rightarrow\infty}\frac{P\{L_n>\log n\}\log n}{EY_n}
    =0,
\]
or equivalently,
\[
    \lim_{n\rightarrow\infty}\frac{EY_n^2}{\mu_n \log n}=0,\quad
    \lim_{n\rightarrow\infty}\frac{P\{L_n>\log n\}\log n}{\mu_n}=0.
\]
The hypothesis (\ref{hyp2}) gives
$$\frac{P\{L_n>\log n\}\log n}{\mu_n}\le
\frac{E(L_n\mathbf{1}_{\{L_n>\log n\}})}
{\mu_n} =o(1)$$
which proves the second desired limit.
For the first limit, for any $\eta\in (0,1)$, write
\begin{align}
EY_n^2&=E[L_n^2\mathbf{1}_{\{L_n\le \eta \log
n\}}]+E(L_n^2\mathbf{1}_{\{\eta\log n<L_n\le\log
n\}})\notag\\ &\le \eta \mu_n \log n+E(L_n\mathbf{1}_
{\{L_n>\eta\log n\}})\log n  \notag \\
&\le (\eta +o_\eta (1))\mu_n \log n
\notag
\end{align}
where we have used (\ref{hyp2}) again to obtain the last inequality.
Thus
$$\frac{EY^2_n}{\mu_n \log n}\le  \eta  + o_\eta(1). $$
Letting $n$ tend to infinity
and then $\eta$ tend to 0 shows that the left-hand side tends to $0$
as desired. \hfill$\Box$\vskip2mm

The next lemma deals with condition
(\ref{E:maint301}) appearing in  Theorem \ref{T:main3} and plays a role in
the proof of Theorem \ref{T:mainc}(2).

\begin{lemma}\label{equiv}
For $n\ge 1$, let $a_n,b_n>0$ and $X_n$ be a non-negative random
variable. According to the sequence $(a_n)_1^\infty$ and $c>0$,
set $Y_n=X_n\mathbf{1}_{\{X_n\le ca_n\}}$ and
$Z_n=Y_n+ca_n\mathbf{1}_{\{X_n>ca_n\}}$. Consider the following
conditions.
\begin{equation}\label{eq1}
    a_n=O(b_n),\,\,\lim_{n\rightarrow\infty}
    \frac{b_nEZ_n^2}{a_n^2EY_n}=0,\,\,
    \lim_{n\rightarrow\infty}\frac{b_n}{EY_n}=\infty.
\end{equation}
Then {\em (\ref{eq1})} holds for some $c>0$ if and only if it holds for
any $c>0$.
\end{lemma}

\begin{proof} On direction is obvious. For the other direction,
we assume that (\ref{eq1}) holds for some $c>0$.
The second condition  in (\ref{eq1}) implies
\begin{equation}\label{E:rmaint3}
    P\{X_n>ca_n\}=o\left(\frac{EY_n}{b_n}\right),\quad
    \frac{EY_n^2}{a_n^2}=o\left(\frac{EY_n}{b_n}\right).
\end{equation}

Let $d>0$ and $Y_n'=X_n\mathbf{1}_{\{X_n\le da_n\}}$ and
$Z_n'=Y'_n+da_n\mathbf{1}_{\{X_n>da_n\}}$. Then (\ref{E:rmaint3})
and Chebyshev inequality imply
\begin{align}
    |EY_n'-EY_n|&\le \begin{cases}ca_nP\{Y_n>da_n\}&
    \mbox{if $d<c$}\\
    da_nP\{X_n>ca_n\}&\mbox{if $d>c$}\end{cases}\notag\\
    & =o\left(\frac{a_nEY_n}{b_n}\right)=o(EY_n)\notag,
\end{align}
and
\begin{align}
    &|EZ_n'^2-EZ_n^2|\le |d^2-c^2|a_n^2P\{
    X_n>(d\wedge c)a_n\}\notag\\
    =&|d^2-c^2|a_n^2\left(P\{Y_n>(d\wedge c)a_n\}
    +P\{X_n>ca_n\}\right)\notag\\
    \le &|d^2-c^2|a_n^2\left(\frac{EY_n^2}{(d\wedge c)^2a_n^2}
    +P\{X_n>ca_n\}\right)
    =o\left(\frac{a_n^2EY_n}{b_n}\right).\notag
\end{align}
Hence we have $EY_n'\sim EY_n$ and
$\frac{b_nE(Z_n')^2}{d^2a_n^2EY_n'}\rightarrow 0$.
\end{proof}

\section{Proofs of Theorems \ref{T:mainc1} and \ref{T:mainc}}

In this section we are concerned with the continuous time process
whose distribution at time $t$, $H_{n,t}$, is given by
(\ref{def-cts}), that is
$$H_{n,t}=e^{-t}\sum_{k=0}^\infty
 \frac{t^k}{k!}Q_{n,p_n}^k.$$
Let $X_{n,1},X_{n,2},...$ be a sequence of independent random variables
with probability distribution $p_n$. Let
$\tilde{X}_n$ be an integer valued random variable
whose probability distribution $\tilde{p}_n$ is given by
\begin{equation}\label{eq-xn*}
  \tilde{p}_n(l)=
    P\{\tilde{X}_n=l\}=\begin{cases}e^{-P\{X_{n,1}\ne 1\}}&
    \text{if $l=1$}\\
    e^{-1}\sum_{1}^{\infty}\frac{1}{j!}
    P\left\{\prod_{1}^j X_{n,i}=l\right\}&
    \text{if $l>1$}\end{cases}.
\end{equation}
With this notation , we have
\[
    H_{n,1}=E(Q_{n,\tilde{X}_n})=Q_{n,\tilde{p}_n}
\]
and
\[
    H_{n,k}=E(Q^k_{n,\tilde{X}_n})=Q^k_{n,\tilde{p}_n},\;\;k=1,2,\dots.
\]

Let $h$ be any nonnegative function defined on $[0,\infty)$
satisfying $h(0)=0$. Fubini's Theorem yields
\begin{equation}\label{E:cgm2}
    E(h(\log\tilde{X}_n))=e^{-1}\sum_{j=1}^{\infty}\frac{1}{j!}
    E(h(\bar{X}_{n,j})),
\end{equation}
where $\bar{X}_{n,j}=\log X_{n,1}+\cdots+\log X_{n,j}$.
Thus, if we assume that $\mu_n,\sigma_n<\infty$ and let
$h(t)=t$ (resp. $h(t)=t^2$), we obtain
\[
    E(\log \tilde{X}_n)=\mu_n\quad \text{and}\quad \text{Var}(\log
    \tilde{X}_n)=\sigma_n^2+\mu_n^2.
\]


\noindent{\bf Proof of Theorem \ref{T:mainc1}.} Here, we deal with the case
where, for each $n$, $p_n(m_n)=1$ for some integer $m_n$.
Observe that for any integers $n,M$ and time $t>0$,
\begin{align}
    &\left\|H_{n,t}-U_n\right\|_{\mbox{\tiny TV}}
    \ge H_{n,t}(id)-\frac{1}{n!}\ge
    e^{-t}-\frac{1}{n!}\notag\\
    &\|H_{n,t}-U_{n}\|_{\mbox{\tiny TV}}\le
    e^{-t}\sum_{i=0}^{M}\frac{t^i}{i!}
    +\|Q_{n,p_{n}}^M-U_{n}\|_{\mbox{\tiny TV}},\notag
\end{align}
where $id$ is the identity of $S_n$, that is, represents the deck
in order.

Assume that
\[
    \liminf_{n\rightarrow\infty}\frac{\log n}{\mu_n}<\infty.
\]
Let $M$ be an integer and  $(n_k)_1^\infty$ be an increasing
sequence such that $\sup_{k\ge 1}\frac{2\log n_k}{\mu_{n_k}}<M$.
Let $(t_k)_1^\infty$ be an arbitrary  sequence of positive
numbers. Then, by Theorem \ref{th-1pt} and the observation above,
we have
\[
    \lim_{k\rightarrow\infty}
    \|H_{n_k,t_k}-U_{n_k}\|_{\mbox{\tiny TV}}=0
    \Longleftrightarrow\lim_{k\rightarrow\infty}t_k=\infty.
\]
This means that the subfamily
$\{(S_{n_k},H_{n_k,t},U_{n_k})\}_1^\infty$, and thus $\mathcal{F}$
itself, does not present a total variation cutoff.

Assume now that
$$\lim_{n\rightarrow \infty}\frac{\log n}{\mu_n}=\infty.$$
Then $t_n=\frac{3\log n}{2\mu_n}$ tends to infinity and
thus $t_n\sim \lfloor t_n\rfloor$. Clearly, a $(t_n,\sqrt{t_n})$ cutoff for
$H_{n,t}$ is equivalent to a $(t_n,\sqrt{t_n})$ cutoff for
$Q^k_{n,\tilde{p}_n}$.
We now prove the desired cutoff
by applying Theorem \ref{T:main2} to $Q_{n,\tilde{p}_n}$.
To this end, we need to show that
(\ref{Lind}) holds for $\tilde{X}_n$.
Set $\tilde{\xi}_n=\frac{\log
\tilde{X}_n-\mu_n}{\sqrt{\sigma_n^2+\mu_n^2}}$. Then
(\ref{E:cgm2}) implies
\[
    E\left(\tilde{\xi}_n^2\mathbf{1}_{\left\{\tilde{\xi}_n^2
    >\epsilon\frac{\log
    n}{\mu_n}\right\}}\right)=\sum_{j>\sqrt{\frac{\epsilon\log
    n}{\mu_n}}}^{\infty}\frac{e^{-1}j^2}{(j+1)!}\rightarrow 0\quad
    \text{as $n\rightarrow\infty$},
\]
for any $\epsilon>0$ and
$n\ge m_n^{1/\epsilon}$.
Hence (\ref{Lind}) holds for $\tilde{X}_n$ and, by Theorem
\ref{T:main2}, the family $\{(S_n,Q_{n,\tilde{p}_n},U_n)\}_1^\infty$
presents, as desired, an optimal $(t_n,b_n)$ total variation cutoff
with $b_n=\sqrt{\log n/\mu_n}$.
\hfill$\Box$\vskip2mm

\noindent{\bf Proof of Theorem \ref{T:mainc}(1).}
As in the proof of Theorem \ref{T:mainc1}, the desired cutoff for
the family $\{(S_n,H_{n,t},U_n)\}_1^\infty$
is equivalent to the same cutoff for
$\{(S_n,Q_{n,\tilde{p}_n},U_n)\}_1^\infty$ because cutoff time and
window size tend to infinity. Hence, the desired conclusion will follow
from Theorem \ref{T:main2} if we can show that
$\tilde{X}_n$ at (\ref{eq-xn*}) satisfies
(\ref{Lind}). Set $\tilde{\xi}_n=\frac{\log
\tilde{X}_n-\mu_n}{\sqrt{\sigma_n^2+\mu_n^2}}$. Then
(\ref{E:cgm2}) implies
\begin{equation}\label{E:pmainc}
    E\left(\tilde{\xi}_n^2\mathbf{1}_{\left\{\tilde{\xi}_n^2
    >\epsilon\frac{\log
    n}{\mu_n}\right\}}\right)=e^{-1}\sum_{j=1}^{\infty}\frac{1}{j!}
    E\left(\frac{\left(\bar{X}_{n,j}
    -\mu_n\right)^2}{\sigma_n^2+\mu_n^2}
    \mathbf{1}_{\left\{\frac{\left(\bar{X}_{n,j}
    -\mu_n\right)^2}{\sigma_n^2+\mu_n^2}>\epsilon\frac{\log
    n}{\mu_n}\right\}}\right),
\end{equation}
if $\epsilon\mu_n^{-1}\log n>1$.
Fix $\epsilon,\delta>0$ and let
$M=M(\delta)\in\mathbb{N}$, $N=N(\epsilon,M)\in\mathbb{N}$ such
that $2\sum_{M+1}^{\infty}\frac{j^2}{j!}<\delta$ and
$\sqrt{\frac{\epsilon\log n}{\mu_n}}\ge 2M$ if $n\ge N$. In this
case, (\ref{E:pmainc}) implies that
\begin{equation}\label{E:pmainc1}
    E\left(\tilde{\xi}_n^2\mathbf{1}_{\left\{\tilde{\xi}_n^2
    >\epsilon\frac{\log
    n}{\mu_n}\right\}}\right)\le
    \delta+e^{-1}\sum_{j=1}^{M}\frac{1}{j!}
    E\left(\frac{\left(\bar{X}_{n,j}
    -\mu_n\right)^2}{\sigma_n^2+\mu_n^2}
    \mathbf{1}_{\left\{\frac{\bar{X}_{n,j}
    -\mu_n}{\sqrt{\sigma_n^2+\mu_n^2}}>\sqrt{\frac{\epsilon\log
    n}{\mu_n}}\right\}}\right).
\end{equation}
To bound the expectation in the right hand side, we consider the
following sets. For $1\le i\le j\le M$, let
\begin{align}
    A_{n,i,j}&=\left\{\log
    X_{n,i}>\frac{1}{j}\left(\mu_n+\sqrt{\frac{\epsilon(\sigma_n^2
    +\mu_n^2)\log
    n}{\mu_n}}\right)\right\}\notag\\
    B_{n,i}&=\left\{\frac{(\log X_{n,i}-\mu_n)^2}{\sigma_n^2}
    >\frac{\epsilon\log
    n}{4M^2\mu_n}\right\}\notag.
\end{align}
Then
\begin{equation}\label{eq-subset}
    \left\{\frac{\bar{X}_{n,j}-\mu_n}{\sqrt{\sigma_n^2+\mu_n^2}}
    >\sqrt{\frac{\epsilon\log n}{\mu_n}}\right\}\subset
    \bigcup_{i=1}^j A_{n,i,j}
\end{equation}
and $$A_{n,i,j}\subset
B_{n,i}\quad\mbox{ if $\sqrt{\frac{\epsilon\log n}{\mu_n}}\ge
2M$}.$$ This implies that for $n\ge N,1\le i\le j\le M$,
\begin{align}
    &E\left(\frac{\left(\bar{X}_{n,j}
    -\mu_n\right)^2}{\sigma_n^2+\mu_n^2}
    \mathbf{1}_{A_{n,i,j}}\right)\notag\\
    \le& 2E\left(\frac{(\bar{X}_{n,j}-\log
    X_{n,i})^2}{\sigma_n^2+\mu_n^2}\mathbf{1}_{B_{n,i}}\right)+
    2E\left(\frac{(\log X_{n,i}-\mu_n)^2}
    {\sigma_n^2}\mathbf{1}_{B_{n,i}}\right)\notag\\
    =& \frac{2\left((j-1)\sigma_n^2+(j-1)^2\mu_n^2\right)}
    {\sigma_n^2+\mu_n^2}P\{B_{n,i}\}+2E\left(\xi_n^2
    \mathbf{1}_{\left\{\xi_n^2>\frac{\epsilon\log
    n}{4M^2\mu_n}\right\}}\right)\notag\\\le& 3E\left(\xi_n^2
    \mathbf{1}_{\left\{\xi_n^2>\frac{\epsilon\log
    n}{4M^2\mu_n}\right\}}\right)\quad\text{if $n$ is large}.\notag
\end{align}

Now, using (\ref{eq-subset})  and these estimates  in
(\ref{E:pmainc1}), and applying the hypothesis that $X_n$ satisfies
(\ref{Lind}), we obtain
\[
    \limsup_{n\rightarrow\infty}
    E\left(\tilde{\xi}_n^2\mathbf{1}_{\left\{\tilde{\xi}_n^2
    >\epsilon\frac{\log
    n}{\mu_n}\right\}}\right)\le
    \delta\quad \forall \delta,\epsilon>0.
\]
Hence (\ref{Lind}) holds for $\tilde{X}_n$.
By Theorem \ref{T:main2},
the family $\{(S_n,H_{n,t},U_n)\}_1^\infty$
presents an optimal $\left(\frac{3\log n}{2\mu_n},b_n\right)$ total
variation cutoff, where
$$b_n=\frac{1}{\mu_n}\max\left\{\sqrt{\frac{(\sigma_n^2+\mu_n^2)
\log n}{\mu_n}},1\right\}$$
(note that $b_n$ always tends to infinity). \hfill$\Box$\vskip2mm

\noindent{\bf Proof of Theorem \ref{T:mainc}(2).}
The proof is similar to that of part (1) except that we will use
Theorem \ref{T:main3} instead of Theorem \ref{T:main2}.
Let
$$\tilde{Y}_n=(\log\tilde{X}_n)\mathbf{1}_{\{\log\tilde{X}_n\le a_n
\}},\;\;
\tilde{Z}_n=\tilde{Y}_n+a_n\mathbf{1}_{\{\log\tilde{X}_n>a_n\}}.$$
By (\ref{E:cgm2}), we have
\[
    E\tilde{Y}_n=e^{-1}\sum_{j=1}^{\infty}\frac{1}{j!}
    E\left[\left(\sum_{i=1}^j\log X_{n,i}\right)
    \mathbf{1}_{\{\sum_{1}^j\log X_{n,i}
    \le a_n\}}\right].
\]
It is apparent that $E\tilde{Y}_n\le EY_n$. For $j>0$, we have
\begin{align}
    E\left[\left(\sum_{i=1}^j\log X_{n,i}\right)
    \mathbf{1}_{\{\sum_{1}^j\log X_{n,i}
    \le a_n\}}\right]\ge \sum_{i=1}^j&\bigg\{
    E\left(\log X_{n,i}\mathbf{1}_{\left\{\log X_{n,i}\le
    \frac{a_n}{j}\right\}}\right)\notag\\
    &\times \prod_{\begin{subarray}{c}k=1\\k\ne i\end{subarray}}^j
    P\left(\log X_{n,k}\le \frac{a_n}{j}\right)
    \bigg\}\notag.
\end{align}
By Lemma \ref{equiv} (or Remark \ref{R4}) and (\ref{E:rmaint3}),
we have
\[
    \liminf_{n\rightarrow\infty}E\left[\left(\sum_{i=1}^j
    \log X_{n,i}\right)
   1_{\{\sum_{1}^j\log X_{n,i}\le a_n\}}\right]\bigg/EY_n\ge j.
\]
Hence, for $k>0$
\[
    \liminf_{n\rightarrow\infty}\frac{E\tilde{Y}_n}{EY_n}\ge
    e^{-1}\sum_{j=0}^k\frac{1}{j!}.
\]
Letting $k\rightarrow\infty$ implies $E\tilde{Y}_n\sim EY_n$.

To apply Theorem \ref{T:main3}, it remains to prove that the
second part of (\ref{E:maint301}) holds for $\tilde{Y}_n$ and
$\tilde{Z}_n$, that is,
\[
    E(\tilde{Y}_n^2)=o\left(\frac{a_n^2EY_n}{\log n}\right),\quad
    P\left\{\log\tilde{X}_n>a_n\right\}=o\left(\frac{EY_n}{\log
    n}\right).
\]
Note that, by  the hypothesis that $X_n$ satisfies
(\ref{E:maint301}), we have
\[
    E(Y_n^2)=o\left(\frac{a_n^2EY_n}{\log n}\right),\quad
    P\left\{\log X_n>a_n\right\}=o\left(\frac{EY_n}{\log
    n}\right).
\]
Then (\ref{E:cgm2}), Lemma \ref{equiv} and the above observation
imply
\begin{align}
    E(\tilde{Y}_n^2)&=e^{-1}\sum_{j=1}^{\infty}\frac{1}{j!}
    E\left[\left(\sum_{i=1}^j
    \log X_{n,i}\right)^2\mathbf{1}_{\{\sum_{1}^j
    \log X_{n,i}\le a_n\}}\right]\notag\\
    &\le e^{-1}\sum_{j=1}^{\infty}\frac{1}{j!}E\left(\sum_{i=1}^j
    (\log X_{n,i})\mathbf{1}_{\{\log X_{n,i}\le
    a_n\}}\right)^2\notag\\
    &=EY_n^2+(EY_n)^2\le 2EY_n^2=
    o\left(\frac{a_n^2EY_n}{\log n}\right)\notag,
\end{align}
and
\begin{align}
    P\{\log \tilde{X}_n>a_n\}&=e^{-1}\sum_{j=1}^{\infty}
    \frac{1}{j!}P\left\{\sum_{i=1}^{j}\log X_{n,i}>a_n\right\}
    \notag\\
    &\le e^{-1}\sum_{j=1}^{\infty}\frac{1}{(j-1)!}P\left\{\log
    X_n>\frac{a_n}{j}\right\}\notag.
\end{align}

Since, for $j\ge 1$,
\begin{align}
    P\left\{\log X_n>\frac{a_n}{j}\right\}&=P\{\log
    X_n>a_n\}+P\left\{Y_n>\frac{a_n}{j}\right\}\notag\\
    &=P\{\log
    X_n>a_n\}+\frac{j^2EY_n^2}{a_n^2}=j^2\times
    o\left(\frac{EY_n}{\log
    n}\right)\notag,
\end{align}
we have
\[
    P\{\log \tilde{X}_n>a_n\}=o\left(\frac{EY_n}{\log n}\right).
\]

By Theorem \ref{T:main3}, the family
$\{(S_n,Q_{n,\tilde{p}_n},U_n)\}_1^\infty$
presents a total variation cutoff
with critical time $\frac{3\log n}{2EY_n}$.
Hence the same holds for $\{(S_n,H_{n,t},U_n)\}_1^\infty$.
\hfill$\Box$


\vskip2mm












\begin{thebibliography}{9}


\bibitem{A} Aldous, D. {\em
Random walks on finite groups and rapidly mixing Markov chains.}
Seminar on probability, XVII,  243--297, Lecture Notes in Math., 986,
Springer, Berlin, 1983.

\bibitem{ADm}
Aldous, D. and  Diaconis, P. {\em Shuffling cards and stopping
times.}  Amer. Math. Monthly  {\bf 93}, 333--348, 1986.

\bibitem{AD} Aldous, D. and Diaconis, P. {\em Strong uniform times and
finite random walks.} Adv. in Appl. Math. {\bf 8}, 69--97, 1986.


\bibitem{BD92}
        Bayer, D. and Diaconis, P.
\emph{Trailing the Dovetail Shuffle to its
        Lair}. Ann. Appl. Probab. {\bf 2} 294--313, 1992.

\bibitem{BHR}
Bidigare, P., Hanlon, P. and  Rockmore, D. {\em A combinatorial
description of the spectrum for the Tsetlin library and its
generalization to hyperplane arrangements.} Duke Math. J. {\bf
99}, 135--174, 1999.


\bibitem{BD} Brown, K. and  Diaconis, P.
{\em Random walks and hyperplane arrangements.} Ann. Probab. {\bf
26}, 1813--1854, 1998.

\bibitem{C06} Chen, G.-Y. {\em Cutoff
phenomenon for finite Markov chains.} Ph.D. dissertation.

\bibitem{Db} Diaconis, P.
{\em  representations in probability and statistics.} Institute of
Mathematical Statistics Lecture Notes---Monograph Series, {\bf
11}. Institute of Mathematical Statistics, Hayward, CA, 1988.

\bibitem{Dc}
Diaconis, P. {\em  The cutoff phenomenon in finite Markov chains.}
Proc. Nat. Acad. Sci. U.S.A.  {\bf 93-4}, 1659--1664, 1996.

\bibitem{D-Dur}
Diaconis, P. {\em Mathematical developments from the analysis of
riffle shuffling.} In ``Groups, combinatorics \& geometry''
(Durham, 2001), 73--97, World Sci. Publishing, 2003.


\bibitem{DFP} Diaconis, P., Fill, J. and Pitman, J. {\em Analysis
of top to random shuffles.} Combin. Probab. Comput. {\bf 1},
135--155, 1992.

\bibitem{D95}
        Durrett, R.  \emph{Probability: Theory and
        Examples}. Duxbury press, 2nd ed, 1995.

\bibitem{M01dr} Mahajan, S. \emph{Shuffles on Coxeter groups.}
Arxiv, 2001.

\bibitem{SC} Saloff-Coste, L. {\em  Lectures on finite Markov chains.}
Lectures on probability theory and statistics (Saint-Flour Summer School,
1996),
301--413, Lecture Notes in Math., 1665, Springer, Berlin, 1997.

\bibitem{SCfg} Saloff-Coste, L. {\em  Random walks on finite groups.}
In Probability on discrete structures, 263--346, Encyclopaedia
Math. Sci., 110, Springer, Berlin, 2004. (H. Kesten, ed.)

\bibitem{S96}
        Shiryaev, A.N. \emph{Probability}. Springer
        Verlag, 2nd ed, 1996.

\bibitem{SGO}  Stark, D.; Ganesh, A. and  O'Connell, N.
{\em  Information loss in riffle shuffling.} Combin. Probab.
Comput. {\bf 11}, 79--95, 2002.

\bibitem{T73}
        Tanny, S.  \emph{A probabilistic interpretation of the Eulerian
        numbers}. Duke Math. J. {\bf 40} 717--722, 1973.
[Correction {\bf
        41} 689, 1974.]





\end{thebibliography}
\end{document}